\documentclass{amsart}

\usepackage{amsfonts,amssymb,amsmath,amsthm}
\usepackage{bbm}
\usepackage{xypic}
\usepackage{stmaryrd}
\usepackage{setspace}
\xyoption{all}

\usepackage{float}
\restylefloat{figure}

\usepackage{tikz}
\usetikzlibrary{arrows}
\usetikzlibrary{matrix}
\usetikzlibrary{positioning}

\usepackage{hyperref}

\usepackage{svn-multi}
\svnid{$Id: twist.tex 661 2011-11-16 11:20:01Z willdonovan $}

\DeclareFontFamily{OT1}{rsfs}{}
\DeclareFontShape{OT1}{rsfs}{n}{it}{<-> rsfs10}{}
\DeclareMathAlphabet{\curly}{OT1}{rsfs}{n}{it}

\newcommand{\smallonebytwo}[2]{\ensuremath{\left( \begin{smallmatrix} #1 & #2 \end{smallmatrix} \right)}}
\newcommand{\smalltwobyone}[2]{\ensuremath{\left( \begin{smallmatrix} #1 \\ #2 \end{smallmatrix} \right)}}

\newcommand{\twobyone}[2]{\ensuremath{\left( \begin{array}{c} #1 \\ #2 \end{array} \right)}}

\newcommand{\varconds}[2]{\ensuremath{\begin{array}{c|c} #1 & #2 \end{array}}}

\newcommand{\Proof}[1]{\smallskip \noindent {\bf Proof:} { \em {#1 \hfill $\Box$ } } }

\newcommand{\ExactTriangle}[3]{\ensuremath{\xymatrix@C=.6em{ & #3 \ar@{-->}[dl] & \\ #1 \ar[rr] & & #2 \ar[ul] }}}
\newcommand{\ExactTriangleWithMaps}[6]{\ensuremath{\xymatrix@C=.6em{ & #5 \ar@{-->}_{#6}[dl] & \\ #1 \ar^{#2}[rr] & & #3 \ar_{#4}[ul] }}}
\newcommand{\FlatExactTriangle}[3]{\ensuremath{\xymatrix{#1\ar[r] & #2 \ar[r] & #3 \ar@{-->}[r] &}}}
\newcommand{\FlatExactTriangleWithMaps}[6]{\ensuremath{\xymatrix{#1\ar^{#2}[r] & #3 \ar^{#4}[r] & #5 \ar@{-->}^{#6}[r] &}}}

\newcommand{\variableLengthPostnikov}[9]{\ensuremath{\xymatrix@!R@R=-3em@C=-3em{ & #1 \ar[rr] \ar[dr] & & #2 \ar[dr]  & \ldots & \ldots & \ldots & #3 \ar[rr] \ar[dr] & & #4 \ar[dr] & \\ #5 \ar[ur] & & #6 \ar@{-->}[ll] \ar[ur] & & #7 \ar@{-->}[ll]& \ldots  & #8 \ar[ur] & & #9 \ar[ur] \ar@{-->}[ll] & & 0 \ar@{-->}[ll]}}}

\newcommand{\leftEndPostnikov}[9]{\ensuremath{\xymatrix@!R@R=-2em@C=-2em{ & #1 \ar[rr] \ar[dr] & & #2 \ar[dr]  & \ldots & \\ #5 \ar[ur] & & #6 \ar@{-->}[ll] \ar[ur] & & #7 \ar@{-->}[ll]& \ldots }}}

\newcommand{\rightEndPostnikov}[9]{\ensuremath{\xymatrix@!R@R=1em@C=1em{ & \ldots & #3 \ar[rr] \ar[dr] & & #4 \ar[dr] & \\ \ldots & #8 \ar[ur] & & #9 \ar[ur] \ar@{-->}[ll] & & 0 \ar@{-->}[ll]}}}

\newcommand{\verticalPostnikov}[9]{\ensuremath{\xymatrix@!R=0.05em{& #5 \ar[ld] \\ #1 \ar[dd] \ar[rd] & \\ & #6 \ar@{-->}[uu] \ar[ld] \\ #2 \ar[rd] & \\ \vdots & #7 \ar@{-->}[uu] \\ \vdots & \vdots \\ \vdots & #8 \ar[ld] \\ #3 \ar[dd] \ar[rd] & \\& #9 \ar@{-->}[uu] \ar[ld] \\ #4 \ar[rd] & \\& 0 \ar@{-->}[uu] \\}}}

\newcommand\C{\mathbb C}

\renewcommand\O{\mathcal O}
\newcommand\PP{\mathbb P}

\newcommand\R{\mathbb R}
\newcommand\Z{\mathbb Z}

\newcommand\HH{\curly H}

\newcommand\To{\longrightarrow}
\newcommand\From{\longleftarrow}
\newcommand\into{\hookrightarrow}

\newcommand\so{\ \ \Longrightarrow\ }

\newcommand\ip{\lrcorner\,}

\newcommand\take{\backslash}
\newcommand\iso{\simeq}

\newcommand\xyhook{\ar@{^{(}->}}
\newcommand\xyhookup{\ar@{^{(}->}}
\newcommand\xyhookdown{\ar@{_{(}->}}

\newcommand\xybend{\ar@/^/}
\newcommand\xydoublebend{\ar@/^2pc/}
\newcommand\xybendup{\ar@/^/}
\newcommand\xybenddown{\ar@/_/}

\newcommand\xyequals{\ar@{=}}

\newcommand\xybirational{\ar@{-->}}
\newcommand\xyquiverA{\ar@{->>}@/^/}

\newcommand\rk{\operatorname{rk}\,}

\newcommand\id{\operatorname{id}}

\renewcommand\Im{\operatorname{Im}}
\newcommand\Hom{\operatorname{Hom}}
\renewcommand\hom{\curly H\!om}
\newcommand\Ext{\operatorname{Ext}}
\newcommand\RHom{\operatorname{\mathbb{R}Hom}}

\newcommand\RDerived{\mathbb{R}}

\newcommand\LDerived{\mathbb{L}}

\newcommand\Ltimes{\overset{\LDerived}{\otimes}}

\newcommand\Perf{\operatorname{Perf}\,}
\newcommand\Supp{\operatorname{Supp}\,}
\newcommand\Coh{\operatorname{Coh}\,}

\newcommand\End{\operatorname{End}}

\newcommand\pt{\operatorname{pt}}

\newcommand\Sym{\operatorname{Sym}}
\newcommand\Tot{\operatorname{Tot}}
\newcommand\codim{\operatorname{codim}}
\newcommand\Cone{\operatorname{Cone}}

\newcommand\mcA{\mathcal A}
\newcommand\mcB{\mathcal B}
\newcommand\mcC{\mathcal C}
\newcommand\mcD{\mathcal D}
\newcommand\mcE{\mathcal E}
\newcommand\mcF{\mathcal F}
\newcommand\mcG{\mathcal G}

\newcommand\mcN{\mathcal N}
\newcommand\mcT{\mathcal T}

\newcommand\quot{/\kern-.7ex/}

\newcommand\beq[1]{\begin{equation}\label{#1}}
\newcommand\eeq{\end{equation}}
\newcommand\beqa{\begin{eqnarray*}}
\newcommand\eeqa{\end{eqnarray*}}

\makeatletter \@addtoreset{equation}{section} \makeatother

\theoremstyle{definition} 
\newtheorem{definition}[equation]{Definition}

\newtheorem{keyDefinition}{Definition}

\theoremstyle{plain}
\newtheorem{theorem}[equation]{Theorem}
\newtheorem{theoremDefinition}[equation]{Theorem / Definition}
\newtheorem{lemma}[equation]{Lemma}
\newtheorem{corollary}[equation]{Corollary}

\newtheorem{proposition}[equation]{Proposition}
\newtheorem{propositionDefinition}[equation]{Proposition / Definition}
\newtheorem{keyProposition}[keyDefinition]{Proposition}

\newtheorem{keyTheorem}[keyDefinition]{Theorem}

\theoremstyle{remark}
\newtheorem{example}[equation]{Example}
\newtheorem{remark}[equation]{Remark}
\newtheorem{remarks}[equation]{Remarks}

\newtheorem{observation}[equation]{Observation}

\newtheorem{fact}[equation]{Fact}

\newcommand{\defined}[1]{{\em\bf #1}}
\newcommand{\emphasis}[1]{{\em #1}}
\newcommand{\remarkEmphasis}[1]{{\em\bf #1}}
\newcommand{\theoremEmphasis}[1]{{\em\bf #1}}
\newcommand{\step}[1]{{\em\bf Step #1:}}
\newcommand{\stepDescribed}[2]{{\em\bf Step #1: (#2)}}

\newcommand{\BibliographyLocation}{../../Bibliography}

\newcommand\bigTriangleNodes[9]{
\def\factor{#1}
\def\xshift{#2}
\node ({#3}Schurk/k) at (-6 * \factor + \xshift, 0) {#4};
\node (above{#3}Schurk/k) at (-4 * \factor + \xshift, 1) {};
\node ({#3}Schur2/2) at (-4 * \factor + \xshift, 0) {\ldots};
\node ({#3}Schur1/1) at (-2 * \factor + \xshift, 0) {#5};
\node ({#3}Schur0/0) at (0 + \xshift,0) {#6};
\node (above{#3}Schur0/0) at (0 + \xshift, 1) {};
\node ({#3}Schurk/k-1) at (-5 * \factor + \xshift, 1 * \factor) {#7};
\node ({#3}Schur2/1) at (-3 * \factor + \xshift, 1 * \factor) {\ldots};
\node ({#3}Schur1/0) at (-1 * \factor + \xshift, 1 * \factor) {#8};
\node ({#3}Schurk/k-2) at (-4 * \factor + \xshift, 2 * \factor) {\ldots};
\node ({#3}Schur2/0) at (-2 * \factor + \xshift, 2 * \factor) {\ldots};
\node ({#3}Schurk/0) at (-3 * \factor + \xshift, 3 * \factor) {#9};
\draw[gray] ({#3}Schur0/0.mid) -- ({#3}Schurk/k.mid) -- ({#3}Schurk/0.mid) -- cycle;
}

\newcommand\twoBigTriangleNodes[6]{
\def\factor{#1}
\def\xshift{#2}
\node ({#3}Schurk+1/k+1) at (-8 * \factor + \xshift, 0) {#4};
\node ({#3}Schurk+1/k) at (-7 * \factor + \xshift, 1 * \factor) {#5};
\node ({#3}Schurk+1/k-1) at (-6 * \factor + \xshift, 2 * \factor) {\ldots};
\node ({#3}Schurk+1/1) at (-5 * \factor + \xshift, 3 * \factor) {#6};
\draw[gray] ({#3}Schur1/1.mid) -- ({#3}Schurk+1/k+1.mid) -- ({#3}Schurk+1/1.mid) -- cycle;
}

\newcommand\bigWindowNodes[1]{
\bigTriangleNodes{#1}{0}{window}{$\O(d-2)$}{$\O(1)$}{$\O$}{$S^\vee (d-3)$}{$S^\vee$}{$\Sym^{d-2} S^\vee$}
}

\newcommand\Schur[1]{\Sigma^{#1} S^\vee}

\newcommand\bigWindowNodesSchurPowers[1]{
\bigTriangleNodes{#1}{0}{window}{$\Schur{d-2,d-2}$}{$\Schur{1,1}$}{$\Schur{0,0}$}{$\Schur{d-2,d-3}$}{$\Schur{1,0}$}{$\Schur{d-2,0}$}
}

\newcommand\bigWindowNodesResults[1]{
\bigTriangleNodes{#1}{0}{results}{$0$}{$0$}{$\O$}{$0$}{$l^\vee$}{$l^{\vee (d-2)}$}
}

\newcommand\twoBigWindowNodes[1]{
\bigWindowNodes{#1}
\twoBigTriangleNodes{#1}{0}{window}{$\O(d-1)$}{$S^\vee (d-2)$}{$\Sym^{d-2} S^\vee (1)$}
}

\newcommand\twoBigWindowNodesArrows[1]{
\bigWindowNodes{#1}
\twoBigTriangleNodes{#1}{0}{window}{$\O(d-1)$}{$S^\vee(d-2)$}{$\Sym^{d-2}S^\vee(1)$}
\draw[thick,->] ({window}Schurk/k) -- ({window}Schurk/k-1) -- ({window}Schurk/k-2) -- ({window}Schurk/0);
\draw[thick,->] ({window}Schurk+1/k+1) -- ({window}Schurk+1/k) -- ({window}Schurk+1/k-1) -- ({window}Schurk+1/1);
\draw[thick,->] ({window}Schur1/1) -- ({window}Schur1/0);
}

\newcommand\twoBigWindowNodesResults[1]{
\bigWindowNodesResults{#1}
\twoBigTriangleNodes{#1}{0}{results}{$l^{\vee (d-2)}$}{$l^{\vee (d-1)}$}{$\O$}

}

\newcommand\twoBigWindowLabelling{
\node[rectangle,draw] (W) at (0,2 * \factor) {$\mcT$};
\draw [very thick,->] (W) -- ({window}Schur1/0);
\node[rectangle,draw] (W') at (-8 * \factor + \xshift,2 * \factor) {$\mcT \otimes \O(1)$};
\draw [very thick,->] (W') -- ({window}Schurk+1/k);
}

\newcommand\inductionStrategyLabelling{
\node[rectangle,draw] (W) at (0,2 * \factor) {$\text{Step 2}$};
\draw [very thick,->] (W) -- ({window}Schur2/0);
\node[rectangle,draw] (W') at (-8 * \factor + \xshift,2 * \factor) {$\text{Step 3}$};
\draw [very thick,->] (W') -- ({window}Schurk+1/k-1);
\node[rectangle,draw] (W'') at (-4 * \factor + \xshift,-1 * \factor) {$\text{Step 1}$};
\draw [very thick,->] (W'') -- ({window}Schur2/2);
}

\newcommand{\Grassmannian}{\mathbbm{Gr}}

\newcommand\cotangentBundle[1]{T^{\vee}{#1}}

\newcommand{\XtiltingObject}{\mathcal{T}}
\newcommand{\baseTiltingObject}{\mathcal{T}_0}
\newcommand{\GtiltingObject}{\mathcal{T}_\Grassmannian}
\newcommand{\twistBase}{X_0}
\newcommand{\maxStratum}{B}
\newcommand{\maxStratumName}{{big stratum}}

\newcommand{\hyperplane}{H}
\newcommand{\leftAdjoint}{L}
\newcommand{\rightAdjoint}{R}
\newcommand{\compose}{}
\newcommand{\composition}{j}

\newcommand{\HomFromV}{A}
\newcommand{\sphericalObject}{\mcE}

\newcommand{\spanningObject}{\mcA}
\newcommand{\spanningPreimageObject}{\mcB}

\newcommand{\littleBundleProjection}{p}
\newcommand{\bigBundleProjection}{p}
\newcommand{\lineBundleProjection}{\hat{q}}

\newcommand{\shiftAbbreviation}{s}
\newcommand{\differential}{\partial}
\newcommand{\unit}{\eta}
\newcommand{\counit}{\epsilon}

\newcommand{\singularBase}{\End^{<r}(V)}
\newcommand{\singularBaseProjection}{\pi_{\operatorname{sing}}}

\newcommand{\fibrePoint}{v}

\newcommand{\CYdimension}{n}
\newcommand{\spanningSetBase}{\Omega'}
\newcommand{\spanningSet}{\Omega}

\begin{document}
\title[\tiny Grassmannian twists on the derived category via spherical functors]{Grassmannian twists on the derived category via spherical functors}
\author{ Will Donovan }
\date{\today}

\begin{abstract}
We construct new examples of derived autoequivalences for a family of higher-dimensional Calabi-Yau varieties. Specifically, we take the total spaces of certain natural vector bundles over Grassmannians $\Grassmannian(r,d)$ of $r$-planes in a $d$-dimensional vector space, and define endofunctors of the bounded derived categories of coherent sheaves associated to these varieties. In the case $r=2$ we show that these are autoequivalences using the theory of spherical functors. Our autoequivalences naturally generalize the Seidel-Thomas spherical twist for analogous bundles over projective spaces.
\end{abstract}

\maketitle
\tableofcontents

\svnid{$Id: introduction.tex 653 2011-11-11 12:21:38Z willdonovan $}

\section{Introduction}
\subsection{Background and construction}

Beginning with the spherical twist of Seidel-Thomas \cite{Seidel:2000} a number of `twist' autoequivalences have now been studied, including \cite{Horja:2001,Huybrechts:2005wn, Anonymous:2008ti, Cautis:2010vi, Addington:2011tla}. These are endofunctors of the bounded derived categories $D^b(X)$ of certain varieties $X$, thought of as mirror to symplectic monodromies \cite{Thomas:2010tg} with Seidel's symplectic twist being the prototypical example \cite{Seidel:1999}. In the simplest case, this symplectic twist is just the topological Dehn twist for an embedded curve in a Riemann surface, hence the name.

In this paper we construct new examples of twist autoequivalences, using the technology of \emphasis{spherical functors} \cite{Anno:2007wo}. We work with a Grassmannian $\Grassmannian = \Grassmannian(r,V)$ of $r$-dimensional subspaces of a $d$-dimensional vector space $V$. Let $S$ denote the tautological subspace bundle and consider the bundle $\Hom(V,S)$, where $V$ denotes a constant bundle:

\begin{keyDefinition} {\bf [Definition \ref{definition.maximal_stratum}]} The total space $X := \Tot(\Hom(V,S))$ is stratified by the rank of the tautological map $V \To \bigBundleProjection^* S$, where $\bigBundleProjection$ denotes the projection $\bigBundleProjection : X \to \Grassmannian$. The \defined{\maxStratumName}, denoted $\maxStratum$, is the locus where the rank of the map is not full. 
\end{keyDefinition}

We will exhibit an autoequivalence of the derived category $D^b(X)$ which we think of as a twist around the big stratum $\maxStratum$.

\begin{remark} $X$ is Calabi-Yau (Section \ref{section.Calabi-Yau}).\end{remark}

\begin{remark}\label{remark.spherical_case} In the case $r=1$, we have $X = \Tot ( V^\vee \otimes \O_{\PP V}(-1) )$ and $\maxStratum = \PP V$, the zero section. We have an autoequivalence of $D^b(X)$ given by the spherical twist around a spherical object, namely the sheaf $\O_\maxStratum$ (see Section \ref{section.spherical_twist} for more details). Our work here generalises this. \end{remark}

We write a point of $X$ as a pair $(S,A)$ with $S \in \Grassmannian$ and $A \in \Hom(V,S)$. Composing $A$ with the inclusion of $S$ we obtain a map $$X \To \End^{\leq r}(V),$$ where $\End^{\leq r}(V)$ denotes the space of endomorphisms of $V$ with rank at most $r$. This map is a resolution of singularities. Restricting it to $\maxStratum$ gives a map $$\singularBaseProjection: \maxStratum \To \singularBase$$ to the space of matrices $\singularBase$ whose rank is less than $r$. This space is singular, and $\singularBaseProjection$ is not flat. We show how we may resolve the spaces $\maxStratum$ and $\singularBase$ and `flatten' the map $\singularBaseProjection$ by a commutative diagram as follows: \beq{diagram.flattening}\xymatrix{\hat{\maxStratum} \ar[r]^f \ar[d]_\pi & \maxStratum \xyhook[r]^i \ar[d]^{\singularBaseProjection} & X \\ \twistBase \ar[r] & \singularBase}\eeq
 
Here $\hat{\maxStratum}$ and $\twistBase$ are smooth, and $\pi$ is a $\PP^{d - r}$-bundle. In particular $\pi$ is flat, which will be crucial in our proof. We then have:

\begin{keyProposition}\label{proposition.Grassmannian_twist} {\bf [Propositions \ref{proposition.F_well-defined}, \ref{proposition.twist_well-defined}]} The functor $$ F := i_* \RDerived f_* \pi^* : D^b(\twistBase) \To D^b(X)$$ is well-defined with a right adjoint $\rightAdjoint$, and there exists a twist functor $T_F$ such that $$T_F\mcA \iso \{ F\rightAdjoint\mcA \overset{\counit_\mcA}{\To} \mcA \}$$ where $\counit$ is the counit of the adjunction, and the braces denote taking the mapping cone.
\end{keyProposition}

We then show:

\begin{keyTheorem}\label{keyTheorem.autoequivalence} {\bf [Theorem \ref{theorem.autoequivalence}]} For $r=2$, the twist functor $T_F$ is an autoequivalence of $D^b(X)$.
\end{keyTheorem}

\begin{remark} In the course of the proof we find that $F$ is a \emphasis{spherical functor}. \end{remark}

\begin{remark} Similar results follow for $r>2$ by an extension of our construction, given in \cite{Donovan:2011ufa}. \end{remark}

\subsection{Relation with other results}

The bundle $\Hom(V,S)$ naturally contains $$\cotangentBundle{\Grassmannian(r,V)} \iso \Hom(V/S,S)$$ as a subbundle, with the inclusion given by composing with the projection $V \to V/S$. Derived autoequivalences of the total space of $\cotangentBundle{\Grassmannian(r,V)}$ have been found by \cite{Cautis:2009vz} in the case that $2r = \dim V$, using methods of $\mathfrak{sl}_2$-categorification: it would be very interesting to establish some connection with their results.

In other directions, we would like to see how our autoequivalences act on the non-commutative desingularisation of $\End^{\leq r}(V)$, recently demonstrated by \cite{Buchweitz:2011ug}. It would also be interesting to understand the mirror to $T_F$: we make some preliminary comments on this in \cite[Section 1]{Donovan:2011vc}.

\subsection{Outline}

We outline the contents of the paper:
\begin{itemize}
\item Section \ref{section.construction} describes the construction of the diagram \eqref{diagram.flattening}, the functor $F$ and the twist $T_F$.
\item Section \ref{section.autoequivalence_criterion} presents a condition, which we refer to as $F$ being \emphasis{Calabi-Yau spherical}, which implies that $T_F$ is an autoequivalence, and allows us to describe its action on the derived category in terms of a \emphasis{spanning set}.
\item Section \ref{section.autoequivalence_property} gives a proof that $F$ is Calabi-Yau spherical in our case.
\item Section \ref{section.twist-properties} describes the action of $T_F$ on $D^b(X)$ and the associated $K$-theory.
\end{itemize}

The Appendices gather technical results needed in our discussion:
\begin{itemize}
\item Appendix \ref{section.existence_of_twist_kernel} explains the Fourier-Mukai techniques needed to define $T_F$.
\item Appendix \ref{section.concepts_of_generation} explains the  tilting generator technology which we use in our proof, and shows in detail that such generators are also \emphasis{split-generators} under suitable assumptions.
\item Appendix \ref{section.tilting} constructs a tilting generator $\XtiltingObject$ on $X$, having recalled the necessary Schur functor formalism.
\item Appendix \ref{section.calculations_on_generator} contains crucial but routine calculations of the action of the relevant functors on $\XtiltingObject$.
\end{itemize}

\subsection{Acknowledgements}

I am grateful to EPSRC for their financial support during my doctoral work, of which this forms a part. I am deeply indebted to many colleagues for useful conversations, suggestions and inspiration, in particular to Nicolas Addington, and to my examiners, Tom Bridgeland and Alessio Corti. Finally I am immensely grateful to my co-supervisor Ed Segal, whose ideas have made much of this work possible, and to my supervisor Richard Thomas for introducing me to this problem and for all his hard work and unflagging support over the years.

\section{Review: Spherical twists}
\label{section.spherical_twist}

In this section we review the spherical twist and spherical functors.

\subsection{Spherical twists (case $r=1$)}
\label{section.spherical_twist_example}

In this simple case we describe the construction of the twist and a proof that it gives an autoequivalence, using a method analogous to our later argument for $r=2$: the main ideas are similar, so we hope it will serve as a guide for the reader.

We have $X = \Tot ( V^\vee \otimes \O_{\PP V}(-1) )$ with maps as follows $$\xymatrix{\PP V \iso \maxStratum\: \xyhook[r]^{\quad i} & X \ar[d]^\littleBundleProjection \\ & \PP V}$$ where $i$ is the inclusion of the zero section. We identify a \emphasis{spherical object} $i_* \O_{\PP V} \in D^b(X)$, and apply the following theorem. For simplicity, we specialise to the case of a Calabi-Yau variety.

\begin{theoremDefinition}\label{theoremDefinition.spherical_object_and_twist} \cite{Seidel:2000} \, Given a variety $X$ of dimension $n$ with $\omega_X \iso \O$, we say that an object $\sphericalObject \in D^b(X)$ is \defined{spherical} if $$\RHom^\bullet_X (\sphericalObject, \sphericalObject) \iso H^\bullet(S^n, \C) \iso \C \oplus \C[-n],$$ where $S^n$ is the topological $n$-sphere. In this case there is an induced autoequivalence $T_\sphericalObject$, the \defined{spherical twist}, given by $$T_\sphericalObject \mcA \iso \{ \sphericalObject \Ltimes \RHom_X(\sphericalObject, \mcA) \To \mcA \}.$$
\end{theoremDefinition}

It follows that to establish an autoequivalence $T_\sphericalObject$ we just have to show:

\begin{lemma}\label{lemma.zero_section_spherical} $\sphericalObject := i_* \O_{\PP V} \in D^b(X)$ is spherical.

\Proof{We need to calculate $\RHom_X(\sphericalObject, \sphericalObject)$. By the adjunction $i_* \dashv i^!$ \cite[Section 3.4]{Huybrechts:2007tf} we have \beqa \RHom_X(i_* \O_{\PP V}, i_* \O_{\PP V}) & \iso & \RHom_{\PP V}(\O_{\PP V}, i^! i_* \O_{\PP V}) \\ & \iso & \RDerived\Gamma_{\PP V} (i^! i_* \O_{\PP V}), \eeqa and therefore we proceed as follows:

\step{1} The zero section $\PP V$ has normal bundle $\mathcal{N} \iso \littleBundleProjection^* (V^\vee \otimes \O_{\PP V}(-1))$. We write this as $V^\vee(-1)$ for brevity, and take the corresponding Koszul resolution: $$\{\wedge^d V(d) \To \wedge^{d-1} V(d-1) \To \ldots \To V(1) \To \underline{\O}\} \overset{\sim}{\To} i_* \O_{\PP V}$$
(The underline indicates the term in the complex which lies in degree $0$.)

\step{2} We then (twisted) restrict the resolution by applying \beqa i^! & =&  \det \mathcal{N}[-\codim i] \otimes i^* \\ & \iso & \wedge^d V^\vee(-d)[-\codim i] \otimes i^*. \eeqa This gives that \beq{equation.restricted_zero_section_resolution} i^! i_* \O_{\PP V} \iso \{ \underline{\O} \To V^\vee(-1) \To \ldots \To \wedge^{d-1} V^\vee(-d+1) \To \wedge^d V^\vee(-d) \}.\eeq

\step{3} We evaluate $\RDerived\Gamma_{\PP V} (i^! i_* \O_{\PP V})$ by taking derived sections of \eqref{equation.restricted_zero_section_resolution}. The middle terms have no cohomology, the left-most term gives just $\C$, and  the right-most term gives $\C[-\dim \PP V]$ by duality. We hence obtain
$$\RDerived\Gamma_{\PP V} (i^! i_* \O_{\PP V}) = \C \oplus \C[-\dim \PP V-\codim i] \iso H^\bullet(S^{\dim X}, \C),$$ as required.
}
\end{lemma}

\begin{remark}\label{remark.spherical_object_proof_plan}
Observe that we:
\begin{enumerate}
\item resolve the spherical object;
\item (twisted) restrict the resolution to the twisting locus $ \maxStratum \iso \PP V$;
\item take derived sections and find that
\subitem the middle terms vanish,
\subitem one required piece comes from sections,
\subitem and the other from higher cohomology by duality.
\end{enumerate}
We will follow a similar plan in our argument: all these steps are reflected, albeit in more complicated ways. (Specifically, we work relative to the base $\twistBase$, and the vanishing becomes more subtle, see Section \ref{section.applying_adjoint}.)
\end{remark}

\subsection{Spherical functors}
\label{section.spherical_functors}

We give a categorical reformulation of the above twist which encompasses our Grassmannian twist $T_F$.

\begin{theoremDefinition}\label{theoremDefinition.spherical_functor_and_twist}  \cite{Anno:2007wo} 
An exact functor $F:\mcD_0 \To \mcD$ between triangulated categories is \defined{spherical} if \begin{enumerate}\item the \defined{cotwist} $$C_F \iso \{ \id \overset{\unit}{\To} RF \}$$ is an autoequivalence of $\mcD_0$, and \item the natural transformation $R \To C_F L$ induced by $R \overset{R\unit}{\To} RFL$ is an isomorphism of functors. \end{enumerate}
In this case there is an induced autoequivalence $T_F$ of $\mcD$, also known as a \defined{spherical twist}, given by $$T_F \iso \{ FR \overset{\counit}{\To} \id \}.$$\end{theoremDefinition}

\begin{remark}For brevity, we leave implicit the requirements that adjoints exist, and that twist and cotwist are well-defined. See \cite{Anno:2007wo} for a full formulation.\end{remark}

\begin{remark}This reduces to the previous theorem by taking $\mcD_0=D^b(\pt)$ and taking $F$ such that $F:\O_{\pt} \longmapsto \sphericalObject$ and $R=\RHom_X(\sphericalObject, -)$ \cite[Section 3, Example 1]{Anno:2007wo}.\end{remark}

\begin{remark}Observe that if our twisting locus $B$ were smooth, and the resolution map $f$ could therefore be taken as the identity, then our twist functor $T_F$ would reduce to a family spherical twist \cite{Horja:2001} \cite[Section 8.4]{Huybrechts:2007tf}. Unfortunately the presence of the map $f$ means that the proof method for the family spherical twist does not transfer to our case. See \cite[Section 2.4]{Donovan:2011vc} for discussion.
\end{remark}


\section{Grassmannian twist construction}
\label{section.construction}
\subsection{Resolving singular strata}
\label{section.resolutions}

In this section we show how to `flatten' the map $\singularBaseProjection$ given in the introduction by constructing the commutative diagram \eqref{diagram.flattening}. First we describe the geometry involved and the resolutions required in our construction. We recall:

\begin{definition} \label{definition.maximal_stratum} The total space $X := \Tot(\Hom(V,S))$ is stratified by the rank of the tautological map $V \To \bigBundleProjection^* S$, where $\bigBundleProjection$ denotes the projection $\bigBundleProjection : X \to \Grassmannian$. The \defined{\maxStratumName}, denoted $\maxStratum$, is the locus where the rank of the map is not full. 
\end{definition}

For $r>1$, $\maxStratum$ is singular, with fibre over a point $S$ of $\Grassmannian$ given by the singular affine cone of homomorphisms \beqa \maxStratum_S & = & \Hom^{< r}(V,S). \eeqa A natural way to resolve this space was suggested to us by \cite{Cautis:2009vz}. Following their notation we write \beqa \maxStratum_S & = & \left\{ \varconds{ \xymatrix{ 0 \ar[r]^r & S \ar[r] & V \ar@/^/[l]^\HomFromV } }{ \rk \HomFromV \leq r-1 } \right\}.\eeqa In \cite{Cautis:2009vz}, inclusions are marked by their codimension: we will omit these when they are clear from context. Now to resolve this space we simply add, for each point, the data of a hyperplane $\hyperplane \subset S$ containing $\Im(\HomFromV)$. This is always possible because $\rk \HomFromV \leq r-1$. We denote the resulting resolution by $$\hat{\maxStratum}_S = \left\{ \xymatrix{ 0 \ar[r]^{r-1} & \hyperplane \ar[r]^1 & S \ar[r] & V \ar@/^/[ll]^\HomFromV } \right\},$$ with the obvious projection map $f_S : \hat{\maxStratum}_S \To \maxStratum_S$. Now we observe:

\begin{lemma} $\hat{\maxStratum}_S$ is smooth.

\Proof{The space of hyperplanes of S, written as $$\{ \xymatrix {0 \ar[r]^{r-1} & \hyperplane \ar[r]^1 & S } \},$$ is just the projective space $\PP^\vee S$, and we may reuse notation and denote its tautological hyperplane bundle by $\hyperplane$. Then $\hat{\maxStratum}_S$ is the total space of the bundle $$\xymatrix{ \Hom(V, \hyperplane) \ar[d] \\ \PP^\vee S. }$$ Everything here is smooth, so we are done.}
\end{lemma}

Observe now that we can perform this construction in a family, by letting $S$ vary as a subspace of a fixed $V$. We then obtain: 

\begin{definition}We have a \defined{resolution} $f: \hat{\maxStratum} \To \maxStratum,$ where we define $$\hat{\maxStratum} :=  \left\{ \xymatrix{ 0 \ar[r]^{r-1} & \hyperplane \ar[r]^1 & S \ar[r] & V \ar@/^/[ll]^\HomFromV } \right\}.$$ The morphism $f$ is the natural one which forgets $H$.
\end{definition}

We can perform a similar construction on $\singularBase$, the space of endomorphisms of rank less than $r$. This gives:

\begin{definition}\label{definition.twist_base_resolution} We have a \defined{resolution} $X_0 \To \singularBase,$ where we define $$\twistBase := \left\{ \xymatrix{0 \ar[r]^{r-1} & \hyperplane \ar[r] & V \ar@/^/[l]^\HomFromV } \right\}.$$ As before the morphism is the one which forgets $H$.
\end{definition}

Putting this all together yields:

\begin{proposition}\label{proposition.resolutions_and_flattening} The resolutions defined above fit into a commutative square: $$\xymatrix{ \hat{\maxStratum} \ar[r]^f \ar[d]_\pi & \maxStratum \ar[d]^{\singularBaseProjection} \\ \twistBase \ar[r] & \singularBase}$$ The map $\pi$ is flat, being the projection map for the bundle $\PP(V/H)$.

\Proof{The maps are the natural forgetful ones, forgetting $\hyperplane$ in the horizontal direction and $S$ in the vertical. The square commutes because forgetting $\hyperplane$ and $S$ in either order gives the same result.

For the last part, we once again reuse the notation $H$ to denote the tautological bundle on $\twistBase$, and we then observe that $\hat{\maxStratum}$ is isomorphic to the total space of the projective bundle $$\xymatrix{ \PP (V/\hyperplane) \ar[d]^\pi \\ \twistBase. } $$ In particular the projection $\pi$ is flat as claimed.
}
\end{proposition}

\subsection{Calabi-Yau property}
\label{section.Calabi-Yau}

We show:

\begin{lemma}\label{lemma.Calabi-Yau} The total space $X$ of our bundle $$\xymatrix{\Hom(V,S) \ar[d]^\bigBundleProjection \\ \Grassmannian(r,V)}$$ is Calabi-Yau.

\Proof{The tangent bundle $T_X$ fits in an exact sequence $$0 \To \bigBundleProjection^* \Hom(V,S) \To T_X \To \bigBundleProjection^* T_\Grassmannian \To 0.$$ We then find \beqa \det T_X & \iso & \det \bigBundleProjection^* T_\Grassmannian \otimes \det \bigBundleProjection^* \Hom(V,S) \\ & \iso & \bigBundleProjection^* ( \det \Hom(S,V/S) \otimes \det \Hom(V,S) ), \eeqa where we use the fact that $T_\Grassmannian \iso \Hom(S,V/S)$. We also have $$0 \To \Hom(V/S,S) \To \Hom(V,S) \To \End(S) \To 0,$$ and so we deduce that $$\det T_X \iso \bigBundleProjection^* \det \End(S).$$ This is trivial because it is self-dual, and so we are done.
}
\end{lemma}

\begin{remark}\label{remark.Sigma-hat-Calabi-Yau} Note that $\twistBase$ is also Calabi-Yau. One way to see this is by setting $r=1$ in the above lemma.
\end{remark}

\subsection{Twist functor definition}
\label{section.twist_construction}

Using Proposition \ref{proposition.resolutions_and_flattening} we now have a diagram of schemes: $$\xymatrix{\hat{\maxStratum} \ar[r]^f \ar[d]_\pi & \maxStratum \xyhook[r]^i \ar[d] & X \\ \twistBase \ar[r] & \singularBase}$$

\begin{definition} We define a functor $F$ as the following composition: $$ F : \qquad \xymatrix{ D^b(\twistBase) \ar[r]^{\pi^*} & D^b(\hat{\maxStratum}) \ar[r]^{\RDerived f_*} & D^b(\maxStratum) \ar[r]^{i_*} & D^b(X). } $$
\end{definition}

For brevity we write $\composition := i \compose f$.

\begin{proposition} $F$ is well-defined, and has a right adjoint given by $$F \dashv \rightAdjoint := \RDerived\pi_* \composition^!.$$
\label{proposition.F_well-defined}

\Proof{As $\pi$ is flat, and $i$ the inclusion of a closed (albeit singular) subscheme $\maxStratum$, we have that $\pi^*$ and $i_*$ are exact functors, and do not have to be derived. Finally the derived functor $\RDerived f_*$ preserves the bounded derived category because $f$ is a proper morphism of noetherian schemes, see \cite[Theorem 3.23 and discussion following]{Huybrechts:2007tf}.

For the existence of the adjoint we observe that $F \iso \RDerived\composition_* \pi^*$ and use the adjunctions \beqa \pi^* & \dashv & \RDerived\pi_*, \\ \RDerived\composition_* & \dashv & \composition^! .\eeqa The second of these is Grothendieck duality \cite[Corollary 3.35]{Huybrechts:2007tf}, which applies as $\composition = i \compose f$ is a composition of proper morphisms, and hence proper \cite[Corollary II.4.8b]{Hartshorne:1977}. Composing the adjoints we obtain an adjoint for $F$.
}
\end{proposition}

We can now properly define our twist:

\begin{propositionDefinition} \label{proposition.twist_well-defined} The \defined{twist} $T_F : D^b(X) \To D^b(X)$ and the \defined{cotwist} $C_F : D^b(\twistBase) \To D^b(\twistBase)$ can be defined as functors of Fourier-Mukai type such that \beqa T_F\mcA & \iso & \{ F\rightAdjoint\mcA \To \mcA \}, \\ C_F \mcB & \iso & \{ \mcB \To \rightAdjoint F \mcB \}. \eeqa The morphisms are induced by the (co)unit of the adjunction $F \dashv \rightAdjoint$.

\Proof{See Appendix \ref{section.existence_of_twist_kernel} for technical details of why suitable kernels exist.}
\end{propositionDefinition}

To end this section we prove a more concrete description of the right adjoint functor $\rightAdjoint$, and the left adjoint $\leftAdjoint$, for use later:

\begin{proposition} \label{proposition.right_adjoint_description} We have $$\rightAdjoint \simeq \RDerived\pi_* (\omega_\pi \otimes -) \LDerived\composition^* [\dim \composition],$$ where we use the relative canonical bundle $\omega_\pi$ given by $\omega_\pi := \omega_{\hat{\maxStratum}} \otimes \pi^* \omega^{-1}_{\twistBase}.$ Furthermore we have a left adjoint $L$ with \beqa \leftAdjoint & \simeq & R [\shiftAbbreviation] \\ & \simeq & \RDerived\pi_* (\omega_\pi \otimes -) \LDerived\composition^* [\dim \pi],\eeqa where $\shiftAbbreviation := \dim \pi - \dim \composition$.

\Proof{We have relative canonical bundles for the morphisms $\pi$ and $\composition$ because the spaces involved are smooth, so we can write \beqa \omega_{\hat{\maxStratum}} & \iso & \omega_\pi \otimes \pi^* \omega_{\twistBase} \\ & \iso & \omega_{\composition} \otimes \composition^*\omega_X. \eeqa Both $X$ and $\twistBase$ are Calabi-Yau (Lemma \ref{lemma.Calabi-Yau}), so we deduce that $\omega_\pi \iso \omega_{\composition}$. Using \cite[Corollary 3.35]{Huybrechts:2007tf} $$\composition^! - \iso \omega_{\composition} [\dim \composition] \otimes \LDerived\composition^* -,$$ the expression for $\rightAdjoint$ is immediate. We now express $\leftAdjoint$ in terms of $\rightAdjoint$ and the Serre functors, denoted $S_{\twistBase}$ and $S_X$, for the categories in question \cite[Remark 1.31]{Huybrechts:2007tf}. By the Calabi-Yau property, the Serre functors are simply shifts: \beqa L & \iso & S_{\twistBase}^{-1} R S_X \\ & \iso & S_{\twistBase}^{-1}  S_X R \\ & \iso & R [\dim X - \dim \twistBase] \\ & = & R [\dim \pi - \dim j]. \eeqa The result follows.}
\end{proposition}

\section{An autoequivalence criterion}
\label{section.autoequivalence_criterion}

We give a condition on a functor $F$ which implies that the corresponding twist $T_F$ is an autoequivalence.

Although we will only apply this to our specific Grassmannian case, we present it in the general triangulated category setting to underscore the formal nature of the proof, and to make the key points more transparent.

\begin{remark} To avoid overwhelming the reader with unnecessary notation, in the following section we write exact triangles of integral functors where we mean triangles of the corresponding Fourier-Mukai kernels \cite{Huybrechts:2007tf}. \end{remark}

\subsection{Calabi-Yau spherical functors}
\label{section.CY_spherical_functors}

\begin{definition}\label{definition.Calabi-Yau_spherical} We say that an integral functor $F : \mcD' \To \mcD$ from an indecomposable, non-trivial triangulated category $\mcD'$ is \defined{Calabi-Yau spherical} if it satisfies:
\begin{enumerate} 
\item Adjoint and twist existence conditions
\begin{itemize}
\item $F$ has integral adjoints $L \dashv F \dashv R$;
\item The adjoint $L$ is full;
\item There exist a twist $T_F$ and a cotwist $C_F$ with adjoints on both sides, such that there exist distinguished triangles $$\ExactTriangleWithMaps{FR}{\counit}{\id,}{}{T_F}{} \qquad \ExactTriangleWithMaps{\id}{\unit}{RF;}{}{C_F}{}$$
\end{itemize}
\item Serre duality conditions
\begin{itemize}
\item $\mcD$ and $\mcD'$ have Serre functors $S$ and $S'$ respectively;
\end{itemize}
\item Compatibility conditions
\begin{itemize}
\item $F$ intertwines $S'$ with an autoequivalence $S^*$ of $\mcD$ (so that $S^* F \iso F S'$);
\item $C_F$ commutes with $S'$;
\end{itemize}
\item Local Calabi-Yau condition
\begin{itemize}
\item $\mcD$ is locally $\CYdimension$-Calabi-Yau with respect to $F$ (in the sense that $SF \iso F[\CYdimension]$) for some $\CYdimension$;
\end{itemize}
\item Sphericity condition
\begin{itemize}
\item There is an isomorphism of functors $R \overset{\sim}{\To} C_F L $ induced by the natural morphism $R \overset{R \unit}{\To} RFL$.
\end{itemize}
\end{enumerate}

\end{definition}

\begin{remark}We will see that a Calabi-Yau spherical functor $F$ is in fact a \remarkEmphasis{spherical functor} in the sense of Anno \cite{Anno:2007wo}. However the above conditions turn out to be easier to check in our situation. This simplifies our work considerably in Section \ref{section.autoequivalence_property}.\end{remark}

\begin{remark}Although the definition is somewhat unwieldy, most of the conditions are immediately satisfied in our case, and should follow very naturally in cases of interest. The final condition is the one we spend almost all of our time proving.
\end{remark}

\begin{remark}The sphericity condition in Definition \ref{definition.Calabi-Yau_spherical} is a direct generalization \cite[Section 8.4]{Huybrechts:2007tf} of Horja's spherical condition: in that case $L$ is always surjective, so the condition translates to $C_F \iso S'[-\CYdimension]$ under our Calabi-Yau assumptions. The intertwinement condition $S^* F \iso FS'$ corresponds to the requirement that Horja's line bundle $\mathcal{L}$ is a restriction from $X$.
\end{remark}

\subsection{Twists of Calabi-Yau spherical functors}

To show that Calabi-Yau spherical functors $F$ give autoequivalences $T_F$, we first show that the cotwist $C_F$ is an autoequivalence.

\begin{remark}\label{remark.spanning_set} Our first step will be to construct a \remarkEmphasis{spanning set} for the triangulated category $\mcD'$. Assuming the existence of a Serre functor, a set $\spanningSetBase \subset \mcD'$ is said to \remarkEmphasis{span} $\mcD'$ if any non-trivial object of $\mcD'$ has a non-trivial $\Hom$ from some element of $\spanningSetBase$. See \cite[Definition 1.47]{Huybrechts:2007tf} for precise definitions.\end{remark}

\begin{lemma}\label{proposition.Calabi-Yau_spherical_functor_cotwist}For $F$ Calabi-Yau spherical, $C_F$ is an autoequivalence.

\Proof{

\stepDescribed{1}{Spanning set for $\mcD'$} We take $$\spanningSetBase := \Im L \cup \ker F.$$ For each $\mcA \in \mcD'$ we seek $\omega \in \spanningSetBase$ such that $\Hom_\mcD(\omega, \mcA) \neq 0$. If $\mcA \in \ker F$ then we simply take $\omega := \mcA$. Otherwise we have $$\Hom_{\mcD'}(LF\mcA, \mcA) \iso \Hom_\mcD(F\mcA, F\mcA) \not\iso 0$$ and we may take $\omega := LF\mcA$. This completes the verification of the spanning set $\spanningSetBase$.

\stepDescribed{2}{Action of $C_F$ on $\spanningSetBase$} We claim that \beqa C_F |_{\ker F} & \iso & [1], \\ C_F |_{\Im L} & \iso & S'[-n]. \eeqa The first follows directly from the definition of $C_F$. For the second we use Serre duality in the form $R \iso S'LS^{-1}$ to note that \beqa C_F L & \iso & R \\ & \iso & S'LS^{-1} \\ & \iso & S'L[-n] \\ & \iso & S'[-n]L \eeqa where we use the left adjoint of the local Calabi-Yau condition, namely $$L S^{-1} \iso L[-n].$$ The claim now follows from the fullness of $L$.

\stepDescribed{3}{Preservation of $\spanningSetBase$ by $C_F$} Using the previous step, it is immediate that $C_F$ takes $\ker F$ to itself, as it simply acts by a shift. Also $C_F$ takes $\Im L$ to itself: this follows by using the left adjoint of the intertwinement assumption, namely $$LS^{*-1} \iso S'^{-1}L,$$ which gives $$C_FL\mcA \iso S'L\mcA[-n] \iso LS^*\mcA [-n].$$

\stepDescribed{4}{Vanishing of $\Hom$s between parts of $\spanningSetBase$} We note that for $\mcB \in \ker F$ we have $$\Hom_{\mcD'}(L\mcA, \mcB) \iso \Hom_\mcD(\mcA, F\mcB) \iso 0.$$ For $\Hom$s in the other direction we use Serre duality to evaluate \beqa \Hom_{\mcD'}(\mcB, L\mcA) & \iso & \Hom_{\mcD'}(S' \mcB, S'LS^{-1} S\mcA) \\ & \iso & \Hom_{\mcD'}(S' \mcB, R S \mcA) \\ & \iso & \Hom_\mcD(FS' \mcB, S \mcA) \\ & \iso & \Hom_\mcD(S^* F \mcB, S\mcA) \\ & \iso & 0, \eeqa where we use our intertwinement assumption $FS' \iso S^*F$.

\stepDescribed{5}{Autoequivalence property} We first note that $C_F$ is integral and therefore exact. Then \cite[Corollary 1.56]{Huybrechts:2007tf} gives the result if
\begin{itemize}
\item $\mcD'$ is indecomposable and non-trivial,
\item $C_F$ has adjoints on both sides and commutes with the Serre functor $S'$,
\item and for all $\omega_i \in \spanningSetBase$ the induced morphism $$\Hom_{\mcD'}(\omega_1, \omega_2) \To \Hom_{\mcD'}(C_F (\omega_1), C_F (\omega_2))$$ is a bijection.
\end{itemize}
The first two conditions follow by assumption, so it remains to check the criterion on $\Hom$s between elements of the spanning set $\spanningSetBase = \Im L \cup \ker F$. The condition holds for $\omega_i \in \Im L$ or $\omega_i \in \ker F$ by Step $2$. The other cases follow from the following $2$ steps, as all $\Hom$s involved vanish. This completes the proof.
}

\end{lemma}

We then deduce:

\begin{proposition}\label{proposition.autoequivance_criterion} For $F$ Calabi-Yau spherical, $T_F$ is an autoequivalence of $\mcD$.

\Proof{We simply note that $F$ is a spherical functor in the sense of Definition \ref{theoremDefinition.spherical_functor_and_twist} by combining the assumptions and the lemma above, and so $T_F$ is an autoequivalence by Theorem \ref{theoremDefinition.spherical_functor_and_twist}.}
\end{proposition}

\subsection{Action on a spanning set}

Here we describe the action of the twist $T_F$ on a spanning set $\spanningSet$ for the triangulated category $\mcD$. We will use this to understand the action of the Grassmannian twist in Section \ref{section.twist-properties}.

\begin{remark}In the definition of a spanning set in Remark \ref{remark.spanning_set} we required every non-trivial object of our category to have a non-trivial $\Hom$ \remarkEmphasis{from} some element of the set. In the presence of a Serre functor we can equivalently require a non-trivial $\Hom$ \remarkEmphasis{to} some element of the set \cite[Exercise 1.48]{Huybrechts:2007tf}. This will be used in the following proposition.\end{remark}

\begin{proposition}\label{proposition.spanning_set} Assume $\mcD$ has a Serre functor $S$, and take a functor $F : \mcD' \To \mcD$ of triangulated categories, with a left adjoint $\leftAdjoint$ as follows: $$\xymatrix{\mcD' \ar@/_/[r]_F^{\perp} & \mcD \ar@/_/[l]_L }$$ Then $$\spanningSet := \Im(F\leftAdjoint) \cup \ker(\leftAdjoint)$$ is a spanning set for $\mcD$. If furthermore $\mcD$ is locally $\CYdimension$-Calabi-Yau with respect to $F$ (in the sense that $SF \iso F[\CYdimension]$) then there are no $\Hom$s between the two parts of $\spanningSet$.

\Proof{To show that $\spanningSet$ spans, take a non-zero object $\spanningObject \in \mcD$. We give a suitable $\omega$ with $\Hom_{\mcD}( \spanningObject, \omega) \not\simeq 0$ in the following cases:

\begin{itemize}
\item {\bf Case $\spanningObject \in \ker(\leftAdjoint)$:} Take $\omega := \spanningObject$ and use $\Hom_{\mcD}(\spanningObject, \spanningObject) \not\simeq 0$.
\item{\bf Case $\spanningObject \notin \ker(\leftAdjoint)$:} We then have $\leftAdjoint \spanningObject \not\simeq 0$ so $$0 \not\simeq \Hom_{\mcD'}(\leftAdjoint \spanningObject, \leftAdjoint \spanningObject) \iso \Hom_{\mcD}(\spanningObject, F\leftAdjoint \spanningObject),$$ and so we may take $\omega := F\leftAdjoint \spanningObject$.
\end{itemize}

This proves that $\spanningSet$ spans. We now show the vanishing of $\Hom$s between the two parts of $\spanningSet$. 

\stepDescribed{1}{backward $\Hom$s} Taking $\spanningObject \in \ker(\leftAdjoint)$ and any $\spanningPreimageObject \in \mcD$ we have $$\Hom_{\mcD}(\spanningObject, F\leftAdjoint\spanningPreimageObject) \iso \Hom_{\mcD'}(\leftAdjoint\spanningObject, \leftAdjoint\spanningPreimageObject) \iso 0,$$ by adjunction.

\stepDescribed{2}{forward $\Hom$s} We similarly observe \beqa \Hom_{\mcD}(F\leftAdjoint\spanningPreimageObject, \spanningObject) & \iso & \Hom_{\mcD}(\spanningObject, SF\leftAdjoint\spanningPreimageObject)^\vee \\ & \iso & \Hom_{\mcD}(\spanningObject, F\leftAdjoint\spanningPreimageObject [\CYdimension])^\vee \\ & \iso & 0, \eeqa where we use the local Calabi-Yau condition and the previous step.}
\end{proposition}

Note in particular that the proposition applies to a Calabi-Yau spherical functor $F$. We now describe the action of the associated twist $T_F$ on $\spanningSet$.

\begin{remark} Here and elsewhere we use the \remarkEmphasis{triangular identities} for the units $\unit$ and counits $\counit$ of our adjunctions. For example for the adjunction $F \dashv R$ we have \beq{equation.triangular}\counit F \circ F \unit = \id_F.\eeq We briefly explain how this arises. The crucial observation is that the functorial adjunction isomorphism $$\psi: \Hom(F- , F-) \overset{\sim}{\To} \Hom(-,RF-)$$ can be explicitly inverted in terms of the counit $\counit$ by $\psi^{-1} := (\counit F) \circ (F -).$ Now \eqref{equation.triangular} follows from the definition of the unit $\unit := \psi(\id_F)$ \cite[Section IV.1]{:SaundersMacLane}. 
\end{remark}

\begin{proposition}\label{proposition.Calabi-Yau_spherical_functor_action}For $F$ Calabi-Yau spherical we have \begin{enumerate} \item $\mcA\in \ker L \so T_F\mcA \iso \mcA,$ and \item $\mcA \in \Im FL \so T_F\mcA \iso S^* \mcA[-\CYdimension+1]$. \end{enumerate}

\Proof{First note that $$L S^{-1} \iso L [-\CYdimension],$$ from the local Calabi-Yau condition by uniqueness of left adjoints.

For the first part, if $\mcA \in \ker L$ then \beqa R\mcA & \iso & S' L S^{-1} \mcA \\ & \iso & S' L \mcA [-\CYdimension] \\ & \iso & 0. \eeqa The result follows by definition of $T_F$.

For the second part, we emulate \cite[Section 8.4]{Huybrechts:2007tf} and observe that by the definitions of $T_F$ and $C_F$ we have a diagram of distinguished triangles:
$$
\xymatrix{
F \ar[r]_{F \unit } \ar@{=}[d] & FRF \ar[r] \ar[d]^{\counit F} & F C_F  \ar@{-->}[r] & \\
F \ar@{=}[r] & F \ar[d] & & \\
& T_F F \ar@{-->}[d] & & \\
& & & \\
}
$$
The commutativity of the top left-hand square follows from a triangular identity for the adjunction $F \dashv R$. Applying the octahedral axiom \cite[Chapter 10]{Kashiwara:2005vw} gives a diagram as follows:
$$
\xymatrix{
F \ar[r] \ar@{=}[d] & FRF \ar[r] \ar[d] & F C_F  \ar@{-->}[r] \ar[d] & \\
F \ar@{=}[r] & F \ar[r] \ar[d] & 0 \ar@{-->}[r] \ar[d] & \\
& T_F F \ar@{-->}[d] \ar^{\sim\quad}[r] & F C_F  [1] \ar@{-->}[d] & \\
& & & \\
}
$$
We then have that \beqa T_F FL & \iso & F C_F L[1] \\ & \iso & F R[1] \qquad\qquad \text{(sphericity condition)} \\ & \iso & F S' L S^{-1} [1]\\ & \iso & S^* F L S^{-1} [1]  \qquad \text{(intertwinement of $S'$ and $S^*$)} \\ & \iso & S^* F L [-\CYdimension+1], \eeqa which yields the result.}

\end{proposition}

\svnid{$Id: autoequivalence-proof.tex 660 2011-11-16 10:37:42Z willdonovan $}

\section{Autoequivalence property for the twist}
\label{section.autoequivalence_property}
\subsection{Orientation}

\begin{remark} From now on we restrict to the case $r=2$, so that the tautological hyperplane bundle $\hyperplane$ on $\hat{\maxStratum}$ is just a line bundle, which we denote by $l$. \end{remark}

We find in this case that:

\begin{observation}The twist base $\twistBase$ given in Definition \ref{definition.twist_base_resolution} is simply the total space of the bundle $$\xymatrix{\Hom(V,l) \ar[d]^p \\ \PP V.}$$ \end{observation}

$C_F$ acts on $D^b(\twistBase)$ by a non-trivial autoequivalence. To understand the action of $C_F$ we use a tilting generator $\baseTiltingObject$ for $D^b(\twistBase)$ given by $$\baseTiltingObject := p^*\left(\O \oplus l^\vee \oplus \ldots \oplus l^{\vee (d-1)}\right).$$ We explain why this is a tilting generator in Section \ref{section.base_tilting_object}.
Our next step is to understand $F$ applied to the summands $p^* l^{\vee k}$. In Section \ref{section.applying_functor_F} we show how to calculate them all at once, using a geometrical method. We use this to deduce the required properties of the cotwist in Section \ref{section.cotwist_on_image}. The proof concludes in Section \ref{section.autoequivalence_proof}.

\subsection{A tilting generator for $X_0$}
\label{section.base_tilting_object}

We give the straightforward proof of the above tilting claim, deferring some other tilting results which we will need until Appendix \ref{section.tilting_object}.

\begin{proposition}\label{proposition.base_tilting_object}$D^b(\twistBase)$ has a tilting generator (Definition \ref{definition.tilting_generator}) given by $\baseTiltingObject := p^* \mcT_\PP$ where $$\mcT_\PP := \O \oplus l^\vee \oplus \ldots \oplus l^{\vee (d-1)} \in \Coh(\PP V).$$

\Proof{Note first that $\twistBase$ is projective over $\singularBase$ by Proposition \ref{proposition.resolutions_and_flattening}, and that $\singularBase$ is a Noetherian affine of finite type, as required in Definition \ref{definition.tilting_generator}.

We now show that $\baseTiltingObject$ is tilting. We have \beqa \RHom_{\twistBase} (p^* \mcT_\PP, p^* \mcT_\PP) & \iso & \RHom_{\PP V}(\mcT_\PP, p_* p^* \mcT_\PP) \\ & \iso & \RHom_{\PP V}(\mcT_\PP, \mcT_\PP \otimes p_* \O_{\twistBase} ) \\ & \iso & \RHom_{\PP V}(\mcT_\PP, \mcT_\PP \otimes \Sym^\bullet \Hom(V,l)^\vee) \\ & \iso & \RHom_{\PP V}\left(\mcT_\PP, \mcT_\PP \otimes \bigoplus_k \left(\Sym^k V \otimes l^{\vee k}\right)\right), \eeqa which splits into terms of the form $$\RHom_{\PP V}(l^{\vee a}, l^{\vee b} \otimes l^{\vee k}) \iso \RDerived\Gamma_{\PP V}(l^{\vee b+k-a}),$$ where $0 \leq a,b \leq d-1$. We note that $b+k-a > -d,$ so Kodaira vanishing gives the result.

We then show that $\baseTiltingObject$ spans $D^b(X)$. By adjunction, we have $$\Hom_X(\baseTiltingObject, -) = \Hom_X(p^* \mcT_\PP, -) \iso \Hom_{\PP V}(\mcT_\PP, p_* -).$$ Now $p$ is affine hence $\bigBundleProjection_*$ is injective, and $\mcT_\PP$ is the Beilinson tilting generator for $\PP V$ \cite[Example 7]{Toda:2008wl}. We deduce that $\baseTiltingObject^{\perp} \iso 0$ and this completes the proof.
}
\end{proposition}

\subsection{Preliminary: pushdowns from resolution $\hat{\maxStratum}$}

We consider the bundle on $\hat{\maxStratum}$ with fibre $l \take \{0\}$, with its natural $\C^*$ action. We have:
$$\xymatrix{\C^* \ar@/^1pc/[r] & l \take \{ 0 \} \ar[d]_\lineBundleProjection & \\
& \hat{\maxStratum} \ar[r]_\composition & X}$$
Now we observe that, by definition of $F$, \begin{eqnarray}\label{equation.F_as_pushforward}
F (l^{\vee k}) & = & \RDerived \composition_* (\pi^* l^{\vee k}) \nonumber \\
& = & \RDerived \composition_* (l^{\vee k}) \nonumber \\
& = & \RDerived \composition_* (\lineBundleProjection_* \O_{l \take \{ 0 \}})_k \nonumber \\
& \iso & \left( \RDerived (\composition\lineBundleProjection)_* \O_{l \take \{ 0 \}} \right)_k.
\end{eqnarray} (The subscript $k$ denotes taking equivariants of weight $k$ for $k \in \Z$, and we note that the bundle $l \take \{0\}$ is a family of affine schemes, so $\lineBundleProjection$ has no higher pushdowns.)

We now define a morphism of schemes $i$ fitting into the following diagram:
\beq{equation.factoring_resolution_pushdown}\xymatrix{l \take \{ 0 \} \ar[d]_\lineBundleProjection \xyhook[r]_i & S \take \{ 0 \} \ar[d]^q \\
\hat{\maxStratum} \ar[r]_\composition & X}\eeq (We write $S \take \{0\}$ for the total space of the tautological bundle $S$, with the zero section removed.)

\begin{definition} The morphism $i$ is defined affine locally (we omit an explicit presentation) so that it maps a closed point $(x,\fibrePoint)$ of the bundle $l \take \{ 0 \}$ which we write as $$ (x,\fibrePoint) = \left( \xymatrix{ 0 \ar[r]^1 & l \ar[r]^1 & S \ar[r] & V \ar@/^/[ll]^\HomFromV }, \qquad 0 \neq \fibrePoint \in l \right),$$ to a closed point $i(x,\fibrePoint)$ of $S \take \{ 0 \}$ given by $$i(x,\fibrePoint) := \left( \xymatrix{ 0 \ar[rr]^2 & & S \ar[r] & V \ar@/^/[l]^{\iota_l \circ\HomFromV }}, \qquad 0 \neq \iota_l (\fibrePoint) \in S \right),$$ where $\iota_l$ denotes the inclusion $l \into S$.
\end{definition}

\begin{lemma}\label{lemma.cut_out_by_section} The map $i$ is a closed embedding, with $\Im i$ cut out scheme\discretionary{-}{-}{-}theoretically by a section $\alpha \in \Gamma(\mcN)$ of the bundle $$\mcN := \Hom(V, \wedge^2 S \{1\}),$$ where $\{1\}$ denotes a shift of weight under the $\C^*$-action.

\Proof{

We first show that $i$ is injective on closed points (the scheme-theoretic result that $i$ is a closed embedding follows by a local calculation, which we omit). If $i$ were not injective, so that $i(x,\fibrePoint) = i(x',\fibrePoint')$ say, then by definition of the map we would have $\iota_l (\fibrePoint) = \iota_{l'} (\fibrePoint')$ for $l \neq l'$, which would imply $\fibrePoint=\fibrePoint'=0$, a contradiction. 

We now define $\alpha$ as the section induced by the following composition of tautological morphisms $$\alpha : V \overset{A}{\To} S \overset{\wedge \fibrePoint}{\To} \wedge^2 S \{1\},$$ between tautological bundles on the bundle $S \take \{ 0 \}$. The map $\alpha$ is zero at a closed point $(x,\fibrePoint)$ precisely when \beq{equation.section_vanishing_criterion} w || \fibrePoint, \quad \forall w \in \Im A. \eeq At such a point $A$ factors through $\langle \fibrePoint \rangle$, hence the point lies in $\Im i$. Conversely, if a point $(x, \fibrePoint)$ is in $\Im i$ then $A$ factors through $l \ni \fibrePoint$ and \eqref{equation.section_vanishing_criterion} holds. Working on the pull-up $(\bigBundleProjection q)^{-1} (U)$ of an open affine $U \subset \Grassmannian$ we see that indeed $\Im i$ is the subscheme of zeroes of $\alpha$.
}
\end{lemma}

Now we observe:

\begin{lemma} The square given in \eqref{equation.factoring_resolution_pushdown} commutes, that is: $$\xymatrix{\ar@{}[dr]|{\Box} l \take \{ 0 \} \ar[d]_\lineBundleProjection \xyhook[r]_i & S \take \{ 0 \} \ar[d]^q \\
\hat{\maxStratum} \ar[r]_\composition & X}$$ In particular, the composite map $\composition \lineBundleProjection$ factors as a closed embedding $i$ followed by a flat projection $q$.

\Proof{The commutativity is clear from the definitions: the horizontal maps forget $l$, and the vertical maps forget $\fibrePoint$.
}
\end{lemma}

\begin{remark} The method used here is similar to that in \cite[Proposition 11.12]{Huybrechts:2007tf}, where the derived pullback via a blow-up map is computed by factoring it into a closed embedding and a flat projection.
\end{remark}

\begin{remark} The embedding $i$ is $\C^*$-equivariant for the natural $\C^*$-action on the bundle $S \take \{0\}$.
\end{remark}

Our calculation now reduces to evaluating \begin{eqnarray}\label{equation.pushforward_is_pushdown_of_structure_sheaf} \RDerived (\composition \lineBundleProjection)_* \O_{l \take \{0\}} & \iso & \RDerived q_* (i_* \O_{l \take \{0\}}) \nonumber \\
& = & \RDerived q_* \O_{\Im i}, \end{eqnarray}
and so to calculate the derived pushdown $\RDerived q_*$ we Koszul resolve $\O_{\Im i}$. To this end, Proposition \ref{proposition.transverse_section} checks that $\alpha$ cuts out the subscheme $\Im i$ \emphasis{transversally}. We begin by observing:

\begin{lemma}$\codim_X B = d-1$.

\Proof{For $S$ fixed we consider the space $\Hom^{\leq \rho}(V,S)$ of homomorphisms with $\rk \leq \rho$. By \cite[Prop 1.1(b)]{Bruns:198wwa} we have that \beqa \dim \Hom^{\leq \rho}(V,S) & = & (\dim V + \dim S)\rho - \rho^2 \\ & = & d+1, \eeqa having set $\rho=1$, and hence the codimension of $\Hom^{\leq \rho}(V,S)$ in $\Hom(V,S)$ is $2d - (d+1) = d-1$. Applying this in a family over the base $\Grassmannian$ gives the result.}

\end{lemma}

\begin{remark}We can see this result explicitly: locally on an open affine $\bigBundleProjection^{-1} (U)$, the subscheme $B$ of $X$ is cut out by $d-1$ independent minors of the $2 \times d$ matrix representing $A \in \Hom(V,S)$.
\end{remark}

We now give a more complete description of $\Im i$:

\begin{lemma}The restriction of $\Im i$ to the fibre over a closed point $x = (S, A) \in X$ is $$\Im i|_{q^{-1} \{ x \}} = \left\{ \begin{array}{cc} S \take \{ 0 \} & x \in \Grassmannian \\ \Im A \take \{ 0 \} & x \in \maxStratum \take \Grassmannian \\ \emptyset & x \in X \take \maxStratum \end{array} \right\}$$ and furthermore $$\dim \Im i = \dim \Grassmannian + d + 2.$$

\Proof{For the first part we use that a closed point $(x, \fibrePoint)$ in $S \take \{0\}$ lies in $\Im i$ precisely when \begin{itemize} \item there exists a line $l \subset S$ such that $\fibrePoint \in l$, and \item $A: V \to S$ factors through the inclusion $\iota_l : l \into S$. \end{itemize}
When $x \in X \take B$ this is impossible, as $A$ is surjective. When $x \in B \take \Grassmannian$ we have $\rk A = 1$ and so we are forced to have non-zero $\fibrePoint \in l = \Im A$. Finally when $x \in \Grassmannian$ we can take any non-zero $\fibrePoint \in l \subset S$, hence the result.

For the second part we decompose $\Im i$ with respect to the natural stratification of $X$ so that $$\Im i = \Im i|_{q^{-1} (B\take\Grassmannian)} \cup \Im i|_{q^{-1} (\Grassmannian)}.$$ Note that $B \take \Grassmannian$ is a large open subset of $B$ and so, using that the fibre of $\Im i$ over a point $x \in B \take \Grassmannian$ has dimension 1, we have \beqa \dim \Im i|_{q^{-1} (B\take\Grassmannian)} & = & \dim B + 1 \\ & = & \dim X - \codim B + 1 \\ & = & (\dim \Grassmannian + 2d) - (d - 1) + 1  \\ & = & \dim \Grassmannian + d + 2, \eeqa which is the dimension claimed for $\dim \Im i$. To conclude we note that for the other stratum we have $$\dim \Im i|_{q^{-1} (\Grassmannian)} = \dim \Grassmannian + 2 < \dim \Grassmannian + d + 2.$$}\end{lemma}

\begin{proposition}\label{proposition.transverse_section} The subscheme $\Im i$ is cut out \theoremEmphasis{transversally} by the section $\alpha \in \Gamma(\mcN)$ of the bundle $\mcN := \Hom(V, \wedge^2 S \{1\}).$

\Proof{This follows from Lemma \ref{lemma.cut_out_by_section}. We use the previous proposition to check that \beqa \codim_{S \take \{0\}} \Im i & = & \dim \Tot({S \take \{0\}}) - \dim \Im i \\ & = & (\dim \Grassmannian + 2d) + 2 - \dim \Im i \\ & = & d, \eeqa as expected. This suffices by smoothness of $S \take \{0\}$.
}
\end{proposition}

\subsection{Technical digression: convolutions}

To present some of our intermediate results more compactly, we choose to use the language of \emphasis{convolutions} for bounded complexes $(\mcA_\bullet, \differential)$ of objects in the derived category. For the reader who does not want to delve into the details, we note the following facts:

\begin{remarks}
\begin{enumerate}
\item Convolutions generalize the mapping cone: for a two-term complex $(\mcA_1 \overset{\differential}{\To} \mcA_0)$ we just have $$\Cone(\mcA_\bullet, \differential) \iso \{\mcA_1 \overset{\differential}{\To} \mcA_0\}.$$
\item Convolutions are defined using \remarkEmphasis{Postnikov systems} \cite[Section IV.2, Exercise 1]{Gelfand} of exact triangles involving the objects and morphisms of the complex. Details are given below.
\item Given a general complex $(\mcA_\bullet, \differential)$ it is not a priori possible to say that the convolution exists or is unique. However if \beq{equation.convolution_condition} \Hom(\mcA_{k+l+1}, \mcA_k[-l]) \iso 0, \qquad k\geq0, \: l>0, \eeq then a unique convolution exists \cite[Section 3.4]{Cautis:2009vz}. In particular if the $\mcA_k$ are sheaves in $\Coh(X)$, this follows immediately from the vanishing of negative $\Ext$s.
\end{enumerate}
\end{remarks}

We now give a formal definition for the case we will require, following \cite[Section 3.4]{Cautis:2009vz}.

\begin{definition}For a complex $(\mcA_\bullet, \differential)$ of objects in $D^b(X)$ whose non-zero terms are given by $$\left(\mcA_d \To \mcA_{d-1} \To \ldots \To \mcA_1 \To \mcA_0\right)$$ we say that $\mcC_d$ is a \defined{convolution} of the complex if there exists a \defined{Postnikov system} as follows:
$$\verticalPostnikov{\mcA_d}{\mcA_{d-1}}{\mcA_1}{\mcA_0}{\mcC_{d}[-d]}{\mcC_{d-1}[-d+1]}{\mcC_{d-2}[-d+2]}{\mcC_1[-1]}{\mcC_0}$$ The triangles on the left are commutative, and those on the right are exact.
\end{definition}

\begin{definition}If the convolution for the complex $(\mcA_\bullet, \differential)$ is unique up to isomorphism, we write it as $\Cone(\mcA_\bullet, \differential)$. \end{definition}

\begin{remark}As for the mapping cone itself, in general uniqueness of convolutions is up to non-unique isomorphism.\end{remark}

\begin{remark}Notice that $\mcA_0 \overset{\sim}{\From} \mcC_0$. It immediately follows that for the two-term complex $(\mcA_1 \overset{\differential}{\To} \mcA_0)$ we have $$\Cone(\mcA_\bullet, \differential) \iso \{ \mcA_1 \overset{\differential}{\To} \mcA_0 \},$$ as expected.\end{remark}

\subsection{Applying the functor $F$}
\label{section.applying_functor_F}
The lemma in this section is analogous to Remark \ref{remark.spherical_object_proof_plan}, step 1, as we describe resolutions for certain sheaves in $\Im F$.

\begin{remark} We identify $Fl^{\vee k}$ as a \remarkEmphasis{non-unique} convolution of a certain complex. However, after applying $R$ to the corresponding Postnikov system in Lemma \ref{proposition.image_within_window} we obtain a complex whose convolution \remarkEmphasis{is} unique. The advantage of this (perhaps unusual) approach is that we avoid the need to keep track of all the data of the Postnikov systems in the intermediate stages: most of the objects will be killed by $R$.\end{remark}

\begin{lemma} \label{proposition.generalized_Koszul_complexes} $ F l^{\vee k}$ is a convolution of a complex of objects $(\mcE_{k,\bullet}, \differential)$ where $$\mcE_{k,j} = \left\{ \begin{array}{cc} \Sym^{k-j} S^\vee & 0 \leq j \leq k,d \\ \Sym^{j-k-2} S(-1)[-1] & 0, k+2 \leq j \leq d \\ 0 & \text{otherwise} \end{array} \right\} \otimes \wedge^j V(j).$$ Here we define $\O(-1) := \wedge^2 S$.

\Proof{The Koszul resolution for $\O_{\Im i}$ on the total space of the bundle $S\take\{0\}$, justified in Proposition \ref{proposition.transverse_section}, gives an isomorphism $$\left\{ \wedge^d \mcN^\vee \overset{ \alpha}{\To} \wedge^{d-1} \mcN^\vee \overset{ \alpha}{\To} \ldots \overset{ \alpha}{\To} \mcN^\vee \overset{ \alpha}{\To} \underline{\O} \right\} \overset{\sim}{\To} \O_{\Im i},$$ where the differentials in the complex are given by wedging with the section $\alpha$, and the underline denotes the degree $0$ term. We consider now the objects corresponding to successive truncations of this complex, as follows: \beqa
\left\{ \underline{\O} \right\} & =: & \mcC_0 \\
\left\{ \mcN^\vee \To \underline{\O} \right\} & =: & \mcC_1 \\
& \vdots & \\
\left\{ \wedge^{d-1} \mcN^\vee \To \ldots \To \mcN^\vee \To \underline{\O} \right\}  & =: & \mcC_{d-1} \\
\left\{ \wedge^d \mcN^\vee \To \wedge^{d-1} \mcN^\vee \To \ldots \To \mcN^\vee \To \underline{\O} \right\}  & =: & \mcC_d  \\
\eeqa
These form a Postnikov system
\beq{equation.koszul_Postnikov} \verticalPostnikov{\wedge^d \mcN^\vee}{\wedge^{d-1} \mcN^\vee}{\mcN^\vee}{\O}{\mcC_{d}[-d]}{\mcC_{d-1}[-d+1]}{\mcC_{d-2}[-d+2]}{\mcC_1[-1]}{\mcC_0}\eeq
with $\mcC_d \iso \O_{\Im i}$. Now we have that $Fl^{\vee k} \iso \RDerived q_* ( \O_{\Im i} )_k$ by the isomorphisms \eqref{equation.F_as_pushforward} and \eqref{equation.pushforward_is_pushdown_of_structure_sheaf}. Therefore applying the functor $\RDerived q_* ( - )_k$ to the system above and writing
$$\mcE_{k,j} :=  \RDerived q_* \left(\bigwedge{}^j \mcN^\vee\right)_k$$
we find that $Fl^{\vee k}$ is a convolution of the complex $\left(\mcE_{k,\bullet}, \differential \right)$. It only remains to show that $\mcE_{k,j}$ takes the form given above. We have \beqa \mcE_{k,j} & = & \RDerived q_* \left(\bigwedge{}^j \Hom(V, \wedge^2 S\{1\})^\vee \right)_k \\ & \iso & \RDerived q_* \left(\O\{-j\} \otimes \wedge^j V \otimes (\wedge^2 S^\vee)^j\right)_k \\ & \iso & \RDerived q_* \left(\O\{k-j\}\right)_0 \otimes \wedge^j V \otimes (\wedge^2 S^\vee)^j, \eeqa for $0 \leq j \leq d$, and knowing the cohomology of $\PP S$ then gives $$\mcE_{k,j} \iso \left\{ \begin{array}{cc} \Sym^{k-j} S^\vee & k \geq j \\ 0 & k = j-1 \\ \Sym^{-2-(k-j)} S \otimes \wedge^2 S[-1] & k \leq j-2 \end{array} \right\} \otimes \wedge^j V \otimes (\wedge^2 S^\vee)^j.$$
This rearranges to give the required result.
}
\end{lemma}

\begin{remark} The differentials $\differential$ are naturally determined by following the Koszul differentials $\ip \alpha$ through the functor $\RDerived q_* ( - )_k$ and the functorial isomorphisms used in the proof. We will not need to do this in our argument, so we omit an explicit description.
\end{remark}

\begin{remark}At least for $k \geq 0$, the convolutions obtained here are in fact examples of \remarkEmphasis{generalized Koszul complexes} \cite{Bruns:1998ty} associated to the degeneracy locus $\maxStratum$ of the tautological map of bundles $V(1) \to S^\vee$ given by the following composition: $$\xymatrix{V(1) \ar[dr] \ar[r]^A & S(1) \ar[d] \ar[r]^\sim & S^\vee \ar[dl] \\ & X & }$$
\end{remark}

\subsection{Applying the adjoint $R$}
\label{section.applying_adjoint}

As in Remark \ref{remark.spherical_object_proof_plan}, steps 2 and 3, we hope to apply $R$ and find that only the first and last terms of the complexes $(\mcE_{k,\bullet}, \differential)$ survive. This does indeed carry through in Lemma \ref{proposition.image_within_window}, with an important caveat: the vanishing only works within a certain range $0 \leq k \leq d-2$.

First we offer some explanation for this phenomenon. The corresponding sheaves $l^{\vee k}$ do not generate $D^b(\twistBase)$, however they do generate the proper subcategory $\Im L$, as recorded in the proposition below.

\begin{proposition}\label{lemma.generators_for_subcategory} $\Im R = \Im L = \langle l^{\vee k} \rangle_{0 \leq k \leq d-2}$

\Proof{We use some technical results from the Appendices. Proposition \ref{proposition.tilting_object} gives us a tilting generator $\XtiltingObject$ for $D^b(X)$, derived from Kapranov's exceptional collection for $\Grassmannian$ \cite{Kapranov:2009wf}. Its summands are explicitly described in Section \ref{section.explicit_tilting_object} and are given by:
\begin{center} \begin{tikzpicture}
\bigWindowNodes{1.3}
\end{tikzpicture} \end{center}
Appendix \ref{section.concepts_of_generation} then gives that $\XtiltingObject$ split-generates $D^b(X)$. Applying $L$ to the summands of $\XtiltingObject$ using Proposition \ref{proposition.action_of_right_adjoint} we obtain
\begin{center} \begin{tikzpicture}
\bigWindowNodesResults{1.3}
\end{tikzpicture} \end{center}
so that all the terms vanish except for $L \Sym^k S^\vee \iso l^{\vee k}$. We then deduce the result for $\Im L$, and the result for $\Im R$ follows similarly (with a shift).
}
\end{proposition}

We also observe:

\begin{proposition}\label{proposition.adjoints_full} The functors $L$ and $R$ are full.

\Proof{We show in Lemma \ref{lemma.preimages_for_Homs_on_X} that the natural map $$\RHom_X(\Sym^a S^\vee, \Sym^b S^\vee) \To \RHom_{\twistBase}(l^{\vee a} , l^{\vee b})$$ induced by the functor $L$ for $a,b \in \Z$ is \emphasis{surjective}. We deduce that $L$ is full. The case of the functor $R$ is similar, with a shift.
}
\end{proposition}

From our new autoequivalence criterion described in Section \ref{section.autoequivalence_criterion}, we see that it suffices to understand the composition $C_FL$ to conclude that $T_F$ is an autoequivalence. Restricting therefore to the generators of $\Im L$ given in the proposition above we have:

\begin{lemma}\label{proposition.image_within_window} For $0 \leq k \leq d-2$, $$RF l^{\vee k} \iso l^{\vee k} \oplus l^{\vee k}[-\shiftAbbreviation],$$ where as before $\shiftAbbreviation := \dim \pi - \dim \composition$.

\Proof{Applying $R$ to the result of Lemma \ref{proposition.generalized_Koszul_complexes} we have that $RF l^{\vee k}$ is a convolution of a complex $(\mcF_{k,\bullet}, \differential)$  where $\mcF_{k,j} := R \mcE_{k,j}$. We claim that only the $\mcF_{k,d}$ and $\mcF_{k,0}$, corresponding to the left-most and right-most terms of the complexes, are non-zero. Now the convolution is defined using a Postnikov system as follows:
$$\verticalPostnikov{\mcF_{k,d}}{\mcF_{k,d-1}}{\mcF_{k,1}}{\mcF_{k,0}}{\mcC_{k,d}[-d]}{\mcC_{k,d-1}[-d+1]}{\mcC_{k,d-2}[-d+2]}{\mcC_{k,1}[-1]}{\mcC_{k,0}}$$ The $\mcC_{k,j}$ are \emphasis{partial convolutions} and $\mcC_{k,d} \iso RFl^{\vee k}$. Our vanishing assumption gives that most of the vertical right-hand maps are isomorphisms and so we see that $$\xymatrix@!C@C=1.6em{\mcC_{k,0} \ar[r]^\sim & \mcC_{k,1} \ar[r]^{\sim\quad} & \quad \ldots \quad \ar[r]^\sim  & \mcC_{k,d-2} \ar[r]^\sim & \mcC_{k,d-1}}.$$ The uppermost distinguished triangle then reads $$\FlatExactTriangle{\mcC_{k,d}[-d]}{\mcF_{k,d}}{\mcF_{k,0}[-d+1]},$$ which gives $$\FlatExactTriangle{\mcF_{k,0}}{RFl^{\vee k}}{\mcF_{k,d}[d]}.$$

 Now we claim specifically that $$\mcF_{k,j} \iso \left\{ \begin{array}{cc} R(\Sym^k S^\vee) & j=0 \\ R(\Sym^{d-k-2}S^\vee(k+1))[-1] & j=d \\ 0 & \text{otherwise} \end{array} \right\} \otimes \wedge^j V.$$ From Proposition \ref{proposition.action_of_right_adjoint} we have that \beq{equation.right_adjoint_vanishing} R(\Sym^{k-j} S^\vee(j)) \iso 0, \quad 0 < j \leq k \leq d-2.\eeq This suffices to show that $R(\mcE_{k,j}) \iso 0$ for $0 < j \leq k$, so it remains to consider $R(\mcE_{k,j})$ with $k+2 \leq j < d$: as we might expect (see Remark \ref{remark.spherical_object_proof_plan}), the vanishing here is dual to the vanishing \eqref{equation.right_adjoint_vanishing}. In this remaining case we have $$\mcE_{k,j} \iso\Sym^{j-k-2} S(-1)[-1] \otimes \wedge^j V(j).$$ Using Lemma \ref{lemma.exceptional_isomorphism}, which follows this one, we have $$\Sym^{j-k-2} S(j-1) \iso \Sym^{j-k-2} S^\vee (k+1).$$ We then see that $R(\mcE_{k,j}) \iso 0$ by applying \eqref{equation.right_adjoint_vanishing} with $$j'=k+1, \: k'=j-1,$$ and verifying that indeed $0 < j' \leq k' \leq d-2$ (under the assumption $k+2 \leq j < d$). Combining all this vanishing with Lemma \ref{proposition.generalized_Koszul_complexes} gives the claim.

Now from the results in Appendix \ref{section.calculations_on_generator}, we see that \beqa \mcF_{k,0} & \iso & R( \Sym^k S^\vee ) \\ & \iso & l^{\vee k} [\dim i - \dim \pi] \\ & = & l^{\vee k} [-s], \\ \mcF_{k,d} & \iso & R( \Sym^{d-k-2} S^\vee (k+1) )[-1] \\ & \iso & l^{\vee k} [\dim i] [-1] \\ & \iso & l^{\vee k} [-d].\eeqa Consequently $RFl^{\vee k}$ is an extension of the following two objects: \beqa \mcF_{k,0} & \iso & l^{\vee k}[-s], \\ \mcF_{k,d}[d] & \iso & l^{\vee k}. \eeqa There is no non-trivial extension of these sheaves, as each $l^{\vee k}$ is a summand of the tilting bundle $\mcT$. This gives the isomorphism. 
}
\end{lemma}

\begin{lemma}\label{lemma.exceptional_isomorphism} $S \iso  S^\vee(-1)$.

\Proof{As $\rk S = 2$ the natural map $S^\vee \otimes \wedge^2 S \To S$ is an isomorphism, and then $\O(-1) :=\wedge^2 S$ gives the result.}
\end{lemma}

\subsection{The cotwist on the image of $L$}
\label{section.cotwist_on_image}

Finally we can characterize the composition $C_FL$, as follows:

\begin{proposition}\label{proposition.RF_within_window} We have $C_FL \iso L[-s]$, indeed there exists a natural isomorphism of functors \beq{equation.natural_between_adjoints}\phi: R \overset{\sim}{\To} C_F L,\eeq induced by the natural transformation \beq{equation.inducing_natural_between_adjoints}R \overset{R\unit}{\To} RFL.\eeq

\Proof{For $\mcA \in D^b(X)$ the component of the claimed natural isomorphism is given by the morphism which makes the following diagram commute: $$\xymatrix{L\mcA \ar[r]^{\unit L_\mcA \quad} & RFL\mcA \ar[r] & C_FL\mcA \ar@{-->}[r] & \\ & R\mcA \ar[u]^{R \unit_\mcA} \ar@{..>}[ur]_{\phi_\mcA} &}$$
It will suffice to check that this is an isomorphism on the summands of our split-generator $\XtiltingObject$ for $D^b(X)$. As in Proposition \ref{lemma.generators_for_subcategory} we use the summands of $\mcT$ given in Proposition \ref{proposition.tilting_object} and described explictly in Section \ref{section.explicit_tilting_object}. As before the only summands which give non-zero objects after applying $\leftAdjoint$ or $\rightAdjoint$ are the $\Sym^k S^\vee$ for $0 \leq k \leq d-2$. By Lemma \ref{proposition.image_within_window} for these the left-hand part of the diagram then reads as follows: $$\xymatrix{l^{\vee k} \ar[r]^{\rho_l\qquad} & l^{\vee k} \oplus l^{\vee k}[-s] \\ & l^{\vee k}[-s] \ar[u]^{\rho_r}}$$ We determine the $\rho$'s. First observe that $$\Hom(l^{\vee k}, l^{\vee k}[-s]) \iso \Hom(l^{\vee k}[-s], l^{\vee k}) \iso 0,$$ and so $$\rho_l = \twobyone{z_l}{0}, \qquad \rho_r = \twobyone{0}{z_r},$$ for $z_l, z_r \in \End(l^{\vee k})$. Now note that our entire setup is invariant, and in particular the morphisms in question, under the $\C^*$-action given by scaling $A$ (this just scales the fibres on all our bundles). Now with this action $l^{\vee k}$ is exceptional in the sense that $\Hom^{\C^*}(l^{\vee k}, l^{\vee k}) \iso \C$.

We next prove that the morphisms are non-trivial. Firstly $\rho_l$ is a component of $\unit$, so it is necessarily non-trivial (otherwise $\counit F \circ F\unit \neq 1$). Secondly $L\unit$ is non-trivial (otherwise $\counit L \circ L\unit \neq 1$) so $R\unit$ is non-trivial by Proposition \ref{proposition.right_adjoint_description}, which gives the result for $\rho_r$. Consequently, using the scaling automorphisms of the $l^{\vee k}$, we have:
$$\xymatrix{l^{\vee k} \ar[r]^{\smalltwobyone{1}{0}\qquad} & l^{\vee k} \oplus l^{\vee k}[-s] \ar[r]^{\quad\smallonebytwo{0}{1}} & l^{\vee k}[-s] \ar@{-->}[r] & \\ & l^{\vee k}[-s] \ar[u]^{\smalltwobyone{0}{1}} \ar@{..>}[ur]_{\phi_{\Sym^k S^\vee}} &}$$ It then follows immediately that $\phi_{\Sym^k S^\vee}$ is an isomorphism. 
}
\end{proposition}

\subsection{Autoequivalence proof}
\label{section.autoequivalence_proof}

\begin{theorem} $T_F$ is an autoequivalence.
\label{theorem.autoequivalence}

\Proof{We claim that $F : D^b(\twistBase) \to D^b(X)$ is \emphasis{Calabi-Yau spherical} as in Definition \ref{definition.Calabi-Yau_spherical} and apply our autoequivalence criterion, Proposition \ref{proposition.autoequivance_criterion}. The previous proposition gives the sphericity condition $R \overset{\sim}{\To} C_F L$. We explain why the other technical conditions hold:

\begin{itemize}
\item The category $D^b(\twistBase)$ is irreducible by \cite[Proposition 3.10]{Huybrechts:2007tf} because $\twistBase$ is smooth, and in particular normal.
\item The existence of $T_F$ and $C_F$ is covered in Section \ref{section.twist_construction}. They possess adjoints because they are of Fourier-Mukai type, and the left adjoint $L$ is full by Proposition \ref{proposition.adjoints_full}.
\item $D^b(X)$ and $D^b(\twistBase)$ are Calabi-Yau categories by Section \ref{section.Calabi-Yau} and hence have Serre functors \beqa S & = & [\dim X], \\ S' & = & [\dim \twistBase], \eeqa respectively. We therefore take $\CYdimension := \dim X$ and $S^* := [\dim X_0]$. The local Calabi-Yau and compatibility conditions are then immediately satisfied, and $S^*$ is clearly an autoequivalence as required.
\end{itemize}

This shows that $F$ is Calabi-Yau spherical, and completes the proof.
}
\end{theorem}

\begin{remark}It follows immediately from the proof of Proposition \ref{proposition.autoequivance_criterion} that $F$ is a \emphasis{spherical functor} as in \cite{Anno:2007wo} and that $C_F$ is an autoequivalence: we give a full description of $C_F$ in \cite{Donovan:2011ufa}.\end{remark}

\section{Properties of the twist}
\label{section.twist-properties}

\subsection{Action on the spanning set $\Omega$}

We observe that Proposition \ref{proposition.spanning_set} applies to $D^b(X)$ to yield a spanning set $\Omega$. The action of the twist $T_F$ on this set is as follows:

\begin{proposition} We have
\begin{enumerate} \item $\mcA \in \ker L \so T_F\mcA \iso \mcA,$ and \item $\mcA \in \Im FL \so T_F\mcA \iso \mcA[-s+1]$. \end{enumerate}

\Proof{This follows immediately from Proposition \ref{proposition.Calabi-Yau_spherical_functor_action}. For the second part we use that by definition $S^*[-\CYdimension+1] \iso [\dim X_0 - \dim X + 1] = [-s+1]$.
}
\end{proposition}

\subsection{Action on $K$-theory}
\label{section.action_on_K-theory}
The autoequivalence $T_F$ induces an endomorphism of the algebraic $K$-theory $K(X)$, which we write $T^K_F$. See \cite[Section 5.2]{Huybrechts:2007tf} for details.

We show now that the spanning set $\Omega$ induces a decomposition of $K(X)$, and exhibit the action of $T^K_F$ on this decomposition.

\begin{definition}\cite[Section 5.2]{Huybrechts:2007tf} We write $K(X)$ for the \defined{algebraic $K$-theory} of $X$. This is the free abelian group generated by the locally free sheaves $\mcE$ on $X$, modulo the equivalence relation that $\mcE \sim \mcF_1 + \mcF_2$ if $\mcE$ is an extension of the $\mcF_i$.
\end{definition}

\begin{definition}For $\mcE^\bullet \in \Perf(X)$ we write $[\mcE^\bullet]$ for its \defined{$K$-theory class} which is given by $$[\mcE^\bullet] := \sum_i (-1)^i \mcE^i.$$
\end{definition}

\begin{remark} Here $\Perf(X) \subseteq D^b(X)$ denotes the subcategory of perfect complexes, that is those complexes isomorphic to bounded complexes of locally frees. For $X$ smooth these categories coincide.
\end{remark}

\begin{remark} A spanning set need not generate the $K$-theory in general. For instance the sheaves $\{ \O_p \}_{p \in \PP^1}$ span $D^b(\PP^1)$ but all have the same $K$-theory class, and yet $\rk (K(\PP^1)) = 2$.
\end{remark}

\begin{remark} Note that $\Omega$ spans $D^b(X)$ but we do \remarkEmphasis{not} make the stronger claim that $\Omega$ \remarkEmphasis{generates} $D^b(X)$. In the latter case it would follow immedately that $[\Omega]$ spans the $K$-theory.
\end{remark}

\begin{lemma}$\rk [\Im FL] = d-1$.

\Proof{From Proposition \ref{lemma.generators_for_subcategory} we have that $\Im L = \langle l^{\vee k} \rangle_k$ where $k$ varies within $0 \leq k \leq d-2$. This gives $\Im FL =  \langle Fl^{\vee k} \rangle_k.$ It remains to check that the $[Fl^{\vee k}]$ are $\Z$-linearly independent in $K(X)$: using Lemma \ref{proposition.image_within_window} we have that \beqa \Im RFL & = & \langle RFl^{\vee k} \rangle_k \\ & = & \langle l^{\vee k} \oplus l^{\vee k} \rangle_k \\ & = & \langle l^{\vee k} \rangle_k \: \subset \: D^b(\twistBase) \eeqa and so immediately we see that $\rk [\Im RFL] = d-1$, and we deduce the result.
}
\end{lemma}

Finally we can show:

\begin{proposition}We have $K(X) \iso [\Im FL] \oplus [\ker L]$ and according to this decomposition $$T^K_F = \begin{pmatrix} \mathbbm{1} & 0 \\ 0 & -\mathbbm{1} \end{pmatrix}.$$

\Proof{$[\ker L]$ is generated by the summands of $\XtiltingObject$ which vanish under $L$: there are $$\rk K(X) - (d-1)$$ of these. The subcategories $\Im FL$ and $\ker L$ are orthogonal under the Euler pairing on $K(X)$ by Proposition \ref{proposition.spanning_set}, and so the decomposition follows. The action of $T^K_F$ then follows from the proposition above.}
\end{proposition}

\appendix
\section{Existence of the twist kernel}
\label{section.existence_of_twist_kernel}

We prove the following technical lemma regarding our functor $F$:

\begin{lemma} The functor $F$ is of Fourier-Mukai type with kernel $$K := \O_{(\pi \times \composition) \hat{\maxStratum}} \in D^b(\twistBase \times X).$$ The kernel is perfect, and its support is proper over both factors $\twistBase$ and $X$.

\Proof{First note that $\pi^*$ is Fourier-Mukai type with kernel given by the graph of the morphism $\pi$ \cite[Exercise 5.4(ii)]{Huybrechts:2007tf}, that is $$\RDerived(\pi \times \id)_* \O_{\hat{\maxStratum}} \in D^b(\twistBase \times \hat{\maxStratum}).$$ Now by \cite[Exercise 5.12(i)]{Huybrechts:2007tf} the composition $\RDerived j_* \pi^*$ is Fourier-Mukai type with kernel \beqa \RDerived (\id \times j)_* \RDerived (\pi \times \id)_* \O_{\hat{\maxStratum}} & \iso & \RDerived (\pi \times j)_* \O_{\hat{\maxStratum}} \\ & \iso & \O_{(\pi \times j)\hat{\maxStratum}}. \eeqa
The last step follows because $\pi \times j$ is a closed embedding. We check this on closed points: by definition $\pi \times j$ takes points to points as follows:
\beqa \hat{\maxStratum} & \To & \twistBase \times X \\ \left( \xymatrix{ 0 \ar[r]^{r-1} & \hyperplane \ar[r]^1 & S \ar[r] & V \ar@/^/[ll]^\HomFromV } \right) & \longmapsto & \left( \xymatrix{ 0 \ar[r]^{r-1} & \hyperplane \ar[r] & V \ar@/^/[l]^\HomFromV }, \xymatrix{ 0 \ar[r]^{r} & S \ar[r] & V \ar@/^/[l]^\HomFromV }\right) \eeqa
For the second part note the kernel is perfect because it is a sheaf on a smooth space, and consequently has a finite resolution by locally frees. The projection maps $p$ and $q$ from $\Supp(K) = (\pi \times \composition)\hat{\maxStratum}$ are shown below:
$$\xymatrix{ & \hat{\maxStratum} \ar[d]_{\wr} \ar[ld]_\pi \ar[rd]^\composition \\ \twistBase & (\pi \times j)\hat{\maxStratum} \ar[l]^{p\quad} \ar[r]_{\quad q} \xyhookdown[d] & X \\ & \twistBase \times X \ar[ul]^p \ar[ur]_q & } $$
They are proper because $\pi$ and $\composition$ are proper. This is clear because $\pi$ is projective, and hence proper, and $\composition = i \circ f$ is a composition of proper morphisms, hence proper.  
}
\end{lemma}

\begin{remark}The assumptions on the kernel suffice to guarantee that the Fourier-Mukai transform and its adjoints preserve boundedness and coherence. We will refer to \cite{Anno:2010ws} for the proof that the twist $T_F$ exists: these assumptions are the ones used there.
\end{remark} 

We can now define our twist:

\begin{proposition} \label{proposition.twist_well-defined_appendix} The twist $T_F : D^b(X) \To D^b(X)$ and the cotwist $C_F : D^b(\twistBase) \To D^b(\twistBase)$ can be defined as functors of Fourier-Mukai type such that \beqa T_F\mcA & \iso & \{ F\rightAdjoint\mcA \overset{\counit_\mcA}{\To} \mcA \}, \\ C_F \mcB & \iso & \{ \mcB \overset{\unit_\mcB}{\To} \rightAdjoint F \mcB \},\eeqa with the morphisms given by the (co)unit of the adjunction $F \dashv \rightAdjoint$.

\begin{remark} We note that the cone construction is non-functorial, so we cannot simply define $T_F$ as the cone on the counit morphism. Instead we follow the standard procedure of constructing a Fourier-Mukai kernel which yields a functor $T_F$ with the required property, and similarly for $C_F$.
\end{remark}

\Proof{

\stepDescribed{1}{twist} From the lemma we have that $F$ is of Fourier-Mukai type with kernel $K$. To obtain a functor $T_F$ as required, we use \cite[Corollary 3.5]{Anno:2010ws} under the assumptions that
\begin{samepage}
\begin{itemize}
\item $K$ is perfect, and
\item $\Supp(K)$ is proper over $\twistBase$ and $X$.
\end{itemize}
\end{samepage}
This gives us a morphism of kernels $Q \overset{\overline{\counit}}{\To} \O_\triangle$ where $\Phi_Q \mcA \iso FR\mcA$ and the following diagram commutes: $$\xymatrix{ \Phi_Q \mcA \ar[r]^{{\Phi_{\overline{\counit}}}_\mcA} \ar[d]_\wr & \Phi_{\O_\triangle}\mcA \ar[d]^\wr \\ FR\mcA \ar[r]_{\counit_\mcA} & \mcA }$$ Here $\counit$ is the counit. The required conditions hold by the lemma above, so we may define $$T_F := \Phi_{\{Q  \overset{\overline{\counit}}{\To} \O_\triangle\}}.$$

\stepDescribed{2}{cotwist} The result for the cotwist follows from the dual result by taking adjoints. Specifically \cite[Theorem 3.1]{Anno:2010ws} similarly gives us a morphism of kernels $Q' \overset{\overline{\counit}}{\To} \O_\triangle$ where $\Phi_{Q'} \mcA \iso LF\mcA$ and the following diagram commutes: $$\xymatrix{ \Phi_{Q'} \mcA \ar[r]^{{\Phi_{\overline{\counit}}}_\mcA} \ar[d]_\wr & \Phi_{\O_\triangle}\mcA \ar[d]^\wr \\ LF\mcA \ar[r]_{\counit_\mcA} & \mcA }$$ We reuse the notation $\counit$ for the counit morphism. Now we use \cite[Proposition 5.9]{Huybrechts:2007tf} to produce kernels which induce right adjoints of the functors $LF$ and $\id$, which are given by: \beqa LF & \dashv & RF, \\ \id & \dashv & \id. \eeqa Noting the Calabi-Yau condition on $X$, the proposition tells us that these are given by applying the functor $\mathbb{D} := (-)^\vee [\dim X]$ to the kernels. We then have: $$\xymatrix{ \Phi_{\mathbb{D} \O_\triangle} \mcA \ar[r]^{{\Phi_{\mathbb{D} \overline{\counit}}}_\mcA} \ar[d]_\wr & \Phi_{\mathbb{D} Q'} \mcA  \ar[d]^\wr \\ \mcA \ar[r]_{\unit_\mcA} & RF \mcA }$$ 
This commutes because the counit $\counit$ and unit $\unit$ are taken to each other by the adjunction isomorphism for the adjunction $LF \dashv RF$. Finally, observing that $\mathbb{D} \O_\triangle \iso \O_\triangle$ we may define $$C_F := \Phi_{\{\O_\triangle  \overset{\mathbb{D} \overline{\counit}}{\To} \mathbb{D} Q'\}}.$$}
\end{proposition}

\section{Concepts of generation}
\label{section.concepts_of_generation}

We clarify two related concepts of generation for the derived category:

\begin{definition} We say that an object $\mcE$ \defined{split-generates} (or simply \defined{generates}) a triangulated category $\mcD$ if the smallest full triangulated subcategory closed under taking direct summands and containing $\mcE$ is $\mcD$ itself.
\end{definition}

Our goal is to show that an object $\mcE$ split-generates $\mcD = D^b(X)$ for a scheme $X$ if $\mcE$ is a \emphasis{tilting generator} in the sense explained below. To do this, we place appropriate smoothness and finite-dimensionality assumptions on $X$.

It will turn out that the tilting generator condition is easy for us to check in our examples. We explain this in the case of $X_0$ in the following Section \ref{section.base_tilting_object}. The case of $X$ is more elaborate, and is deferred to Appendix \ref{section.tilting}.

\begin{definition}\label{definition.tilting_generator} (cf. \cite[Definition 6]{Toda:2008wl})
We say that a locally free sheaf $\mcE$ on a scheme $X$, where $X$ is projective over a Noetherian affine of finite type, is a \defined{tilting generator} for $D^b(X)$ if 
\begin{enumerate} 
\item $\mcE$ is \defined{tilting} in that it satisfies $\RHom^{>0}_X(\mcE,\mcE) \iso 0$;
\item $\mcE$ is \defined{spanning} in the sense that $0 \iso \mcE^{\perp} \subset D^-(X)$.
\end{enumerate}
\end{definition}

It is standard that:

\begin{proposition}For $\mcE$ a tilting generator as above there exist quasi-inverse equivalences $$ \xymatrix@C=3em{D^b(A) \xybend[rr]^{\quad \Psi} & & D^b(X) \xybend[ll]^{\quad \Phi}} $$ Here $A := \End_X(\mcE)$, we write $D^b(A):=D^b(A\operatorname{-mod})$, and we define \beqa \Psi(-) & := & - \underset{A}{\Ltimes} \mcE, \\ \Phi(-) & := & \RHom_X(\mcE,-).\eeqa 

\Proof{This is \cite[Lemma 8]{Toda:2008wl}.}
\end{proposition}

\subsection{Boundedness for tilting functors}

We now record some boundedness properties of the functors $\Phi$ and $\Psi$:

\begin{proposition}\label{proposition.boundedness_for_Phi} For $\mcF^\bullet \in D^b(X)$ we have: $$\HH^{\geq m}\mcF^\bullet \iso 0 \so \HH^{\geq m+\dim X}\Phi(\mcF^\bullet) \iso 0.$$

\Proof{To evaluate $\Phi(\mcF^\bullet)$ we use $$\Phi(-) := \RHom_X(\mcE, -) \iso \R\Gamma_X \hom_X(\mcE, -).$$ The $\hom$ need not be derived because $\mcE$ is a locally free sheaf. Our assumption $\HH^{\geq m}\mcF^\bullet \iso 0$ implies that $\HH^{\geq m}\hom_X(\mcE, \mcF^\bullet) \iso 0$. Now for any sheaf $\mcG$ we have $\HH^{>\dim X} \R\Gamma_X (\mcG) \iso 0$ by Grothendieck vanishing \cite[Theorem III.2.7]{Hartshorne:1977}, and so the result follows from vanishing in the following spectral sequence for $\R\Gamma_X$ \cite[Equation 2.6]{Huybrechts:2007tf}: $$E^2_{p,q} = \HH^p \R\Gamma_X(\HH^q (-)) \so \HH^{p+q} \R\Gamma_X(-).$$}
\end{proposition}

\begin{proposition}\label{proposition.boundedness_for_Psi} For $M \in A\operatorname{-mod}$ we have $$\HH^{>0}\Psi(M) \iso 0,$$ $$\HH^{<-\dim X}\Psi(M) \iso 0.$$

\Proof{The first vanishing follows directly from the definition of $\Psi$. We show how to deduce the second from Proposition \ref{proposition.boundedness_for_Phi}. Following \cite[Lemma 8]{Toda:2008wl}, we consider the canonical map $$\rho : \tau_{< m} \Psi(M) \To \Psi(M)$$ where $\tau_{< m}$ is a truncation functor. This is defined by $\tau_{<m} := \tau_{\leq m-1}$ where $$(\tau_{\leq n} \mcF^\bullet)^i := \left\{ \begin{array}{ll}\mcF^i & i < n \\ \ker \differential & i = n \\ 0 & i > n \end{array} \right\}$$ as in \cite[Section 1]{Toda:2008wl}. The crucial property of this functor for us is that $$\HH^{\geq m} \tau_{< m} \mcF^\bullet \iso 0,$$ whereas $\HH^i \rho$ is an isomorphism for $i < m$. Now applying $\Phi$ we obtain $$\Phi(\rho) : \Phi(\tau_{< m} \Psi(M)) \To \Phi\Psi(M) \iso M.$$ If we put $m:=-\dim X$ then Proposition \ref{proposition.boundedness_for_Phi} gives that $$\HH^{\geq 0} \Phi(\tau_{< m} \Psi(M)) \iso 0,$$ and then we see that $\Phi(\rho)$ must be zero, as its codomain $M \in D^b(A)$ is a complex concentrated in degree $0$. It follows that $\rho$ itself is zero, as $\Phi$ is an equivalence. Applying $\HH^i$ to $\rho$ for $i < -\dim X$ then allows us to deduce that $\HH^{<-\dim X} \Psi(M)\iso 0$ as required.}
\end{proposition}

\subsection{Consequences of smoothness}

Now assuming furthermore that $X$ is smooth we obtain the following lemma:

\begin{lemma}\label{lemma.higher_ext_vanishing} Given $X$ as above and additionally \theoremEmphasis{smooth}, for $M\in A\operatorname{-mod}$ we have $$\Ext^i_A(M,-) \iso 0$$ for $i \gg 0$, where the placeholder stands for an element of $A\operatorname{-mod}$.

\Proof{We note that \beq{equation.exts_by_smoothness}\Ext^i_A(M,-) \iso \HH^i \RHom_A(M,-) \iso \HH^i \RHom_X(\Psi(M),\Psi(-))\eeq using the equivalence $\Psi$. We want to show that this functor vanishes for $i \gg 0$:

\step{1} We show that for sufficiently large $i$ we have that $$\HH^i \RHom_X(\Psi(M),-)\iso 0,$$ where the placeholder now stands for a coherent sheaf on $X$. Proposition \ref{proposition.boundedness_for_Psi} gives that $\HH^{<-\dim X} \Psi(M) \iso 0$. We now use the spectral sequence \cite[Equation 2.8]{Huybrechts:2007tf} $$E^{p,q}_2 = \HH^p \RHom_X(\HH^{-q} \Psi(M), -) \so \HH^{p+q} \RHom_X(\Psi(M), -).$$ Any coherent sheaf on $X$ has a locally free resolution of length at most $\dim X+1$ by smoothness \cite[Proposition 3.26, and remarks following]{Huybrechts:2007tf}, and it follows that there exists $N$ such that $$\HH^{>N}\RHom_X(-,-)\iso 0,$$ where once again the placeholders stand for coherent sheaves on $X$. (We see this by using the locally free resolutions to evaluate the $\RHom$.) The resulting vanishing in the spectral sequence above suffices to deduce that $$\HH^{>\dim X+N} \RHom_X(\Psi(M),-)\iso 0,$$ as required.
 
\step{2} Now we consider the spectral sequence \cite[Equation 2.7]{Huybrechts:2007tf} $$E^{p,q}_2 = \HH^p \RHom_X(\Psi(M), \HH^q\Psi(-)) \so \HH^{p+q} \RHom_X(\Psi(M), \Psi(-)).$$ We have that $\HH^{>0} \Psi(-) \iso 0$ and so the previous step gives that \eqref{equation.exts_by_smoothness} vanishes for $i>\dim X+N$, and we are done. }
 \end{lemma}

\begin{remark}We briefly indicate how locally free resolutions of length $\dim X+1$ for coherent sheaves $\mcG$ on $X$ are obtained. By \cite[Exercise III.6.8]{Hartshorne:1977} we can construct a locally free resolution $\mcF^\bullet \to \mcG$. We can then truncate this to give an acyclic complex $$0 \To \Im \differential \To \mcF^{-\dim X+1} \To \ldots \To \mcF^0 \To \mcG \To 0.$$ It follows by smoothness of $X$ that $\Im \differential$ is in fact locally free \cite[Proof of Lemma 2.5]{Burban:kWm0LpP3}, so this yields the required resolution.
\end{remark}

We then have:

\begin{proposition}\label{proposition.tilting_generator_criterion} If $\mcE$ is a tilting generator for $D^b(X)$ for $X$ smooth, then $\mcE$ also split-generates $D^b(X)$.

\Proof{Considering $A \in A\text{-mod}$ we have that $$\Phi(A) = A \underset{A}{\Ltimes} \mcE \iso \mcE$$ and we deduce that $\mcE$ split-generates $D^b(X)$ precisely when $A$ split-generates $D^b(A)$. We prove the latter claim as follows:

\step{1} We use the lemma to deduce that every $A$-module $M$ has a \emphasis{finite} projective resolution. For this we first note that the category of $A$-modules has enough projectives, so every $A$-module $M$ has a resolution by projective $A$-modules. Following \cite[Section III.5.9]{Gelfand} we write $\operatorname{pdim} M$ for the largest integer $i$ such that $\Ext^i_A(M,-) \not\iso 0$: this exists because of the smoothness of $X$ by Lemma \ref{lemma.higher_ext_vanishing}. Using \cite[Corollary III.5.12(a)]{Gelfand}, we find that $M$ then has a projective resolution of length $\operatorname{pdim} M + 1$.

\step{2} The previous step can be used to yield finite projective resolutions of more general objects $M^\bullet$ in $D^b(A)$. These are given by bounded complexes of $A$-modules $M^i$. We may resolve each $M^i$ separately to produce a bounded double complex \cite[Section 11.5]{Kashiwara:2005vw} of projective $A$-modules: the total complex of this is then a finite projective resolution of $M^\bullet$.

\step{3} Finally we show that $A$ split-generates $D^b(A)$. Consider then the smallest full triangulated subcategory $\mcC$ of $D^b(A\text{-mod})$ closed under taking direct summands and containing $A$. This contains the free $A$-modules $A^{\oplus i}$ (as these are iterated extensions of $A$), and the projective $A$-modules (as these are direct summands of the frees). It then follows from the previous step that $\mcC = D^b(A\text{-mod})$. This suffices to conclude.
}
\end{proposition}
\section{Tilting generator construction}
\label{section.tilting}
\subsection{Schur functors}

We briefly introduce \emphasis{Schur functors}, as they occur in the description of our tilting generator:

\begin{definition} \defined{(Schur functor for vector spaces)} Take a weight $\lambda$ for $GL(W)$. We write $\Sigma^\lambda W$ for the \defined{Schur power} which is the $GL(W)$-representation with highest weight $\lambda$, or $0$ if such a representation does not exist.
\end{definition}

\begin{remark} Weights for $GL(W)$ correspond to sequences $( \lambda_1, \ldots, \lambda_{\dim(W)} )$ of integers, and are ordered lexicographically. The weights occurring as \remarkEmphasis{highest weights} in $GL(W)$-representations are given by \remarkEmphasis{non-increasing} sequences of integers.
\end{remark}

\begin{example} Take $W := V$, $\dim V = 4$. Then we have for example:
\beqa
\Sigma^{1,0,0,0} V & = & V \\
\Sigma^{1,1,0,0} V & = & \wedge^2 V \\
\Sigma^{1,1,1,0} V & = & \wedge^3 V \\
\Sigma^{1,1,1,1} V & = & \wedge^4 V \\
\Sigma^{k,0,0,0} V & = & \Sym^k V
\eeqa
It follows from the definitions that the representations given have the required highest weight: for their irreducibility we refer to \cite{Fulton:1996tk}.

For a general dominant weight, the description of the Schur functor will be more complex. For an example of a non-dominant weight we have:
\beqa
\Sigma^{0,0,0,1} V & = & 0
\eeqa
\end{example}

In general the rule for multiplying Schur powers is quite elaborate, however we quote:

\begin{fact} \defined{(Pieri formula)} Given a weight of the form $\mu = (1, \ldots, 1, 0, \ldots, 0)$ we have $$\Sigma^\lambda W \otimes \Sigma^\mu W = \bigoplus_{\mu' \in S_d(\mu)} \Sigma^{\lambda + \mu'} W,$$ where $\mu'$ ranges over the orbit of $\mu$ under the natural permutation action of $S_d$, where $d=\dim W$. Weights are added component-wise. \cite[Appendix A, equation A.7]{Fulton:1996tk}
\end{fact}

\begin{definition} {\bf (Schur functor for vector bundles)} Take a vector bundle $E$ with structure group $GL(W)$. We can view $E$ as a principal $GL(W)$-bundle via the frame bundle construction. Given a weight $\lambda$ of $GL(W)$, we defined the \defined{Schur power} $\Sigma^\lambda E$ by $$\Sigma^\lambda E := E \underset{GL(W)}{\otimes} \Sigma^\lambda W.$$
\end{definition}

\begin{example}\label{example.Schur_expansion_for_Grassmannian} Take $E := S^\vee$, the dual of the tautological subspace bundle on the Grassmannian $\Grassmannian(2,V)$. The structure group here is $GL(2)$, and the highest weights are given by pairs $(\lambda_1, \lambda_2)$ with $\lambda_1 \geq \lambda_2$. In this simple case, the Pieri rule shows that the Schur powers decompose into products of $\Sym$s and $\wedge$s and we have: \beqa
\Sigma^{1,0} S^\vee & = & S^\vee \\
\Sigma^{1,1} S^\vee & = & \wedge^2 S^\vee = \O(1) \\
\Sigma^{2,1} S^\vee & = & \Sigma^{1,0} S^\vee \otimes \Sigma^{1,1} S^\vee = S^\vee(1) \\
\Sigma^{2,2} S^\vee & = & \Sigma^{1,1} S^\vee \otimes \Sigma^{1,1} S^\vee = \O(2) \\
\Sigma^{2,0} S^\vee & = & \Sym^2 S^\vee \\
\vdots & & \vdots
\eeqa (We work with the dual $S^\vee$ so that signs match between the left- and right-hand sides under our chosen polarization.) 
\end{example}

\subsection{Construction}
\label{section.tilting_object}

We quote:

\begin{proposition}\label{proposition.kapranov_exceptional_collection}\cite[Section 3]{Kapranov:2009wf} There exists a full strong exceptional collection for $D^b(\Grassmannian(r,V))$ given by suitable Schur powers $$\{ \Sigma^\alpha S^\vee \}_{0 \leq \alpha \leq \alpha_{top} },$$ for $GL(r)$-weights $\alpha$ where $\alpha_{top} := (d-r, \ldots, d-r)$, for $d=\dim V$.
\end{proposition}

By standard arguments this yields a tilting generator for the Grassmannian $\Grassmannian$: $$\GtiltingObject := \bigoplus_{0 \leq \alpha \leq \alpha_{top} } \Sigma^\alpha S^\vee \in D^b(\Grassmannian).$$ We then obtain:

\begin{proposition}\label{proposition.tilting_object} There exists a tilting generator for $D^b(X)$ (Definition \ref{definition.tilting_generator}), constructed by pullback from the base $\Grassmannian$ as follows: $$\XtiltingObject := \bigBundleProjection^* \GtiltingObject = \bigBundleProjection^* \left( \bigoplus_{0 \leq \alpha \leq \alpha_{top} } \Sigma^\alpha S^\vee \right).$$

\Proof{Using a similar approach to \cite[Proposition 4.1]{Bridgeland:2005wi}, we first show that $\XtiltingObject$ is tilting, then demonstrate that it spans the derived category. This is an elaboration of Proposition \ref{proposition.base_tilting_object}. As indicated in the introduction, $X$ is a resolution of the affine singularity $\End^{\leq r}(V)$: in particular it is projective over a Noetherian affine of finite type as required in Definition \ref{definition.tilting_generator}.
\\
\step{1} To show that $\XtiltingObject$ is tilting we require $\RHom^{>0}_X (\XtiltingObject, \XtiltingObject) \iso 0$. First observe that \beqa \RHom_X(\XtiltingObject,\XtiltingObject) & \iso & \RHom_X\left(p^* \GtiltingObject, p^*\GtiltingObject\right) \\ & \iso & \RHom_X\left(\GtiltingObject, p_* p^*\GtiltingObject\right) \\ & \iso & \RDerived\Gamma_X\left(\GtiltingObject^\vee \Ltimes p_* p^*\GtiltingObject\right) \\ & \iso & \RDerived\Gamma_X\left(\GtiltingObject^\vee \Ltimes\GtiltingObject\Ltimes p_* \O_X \right) \\ & \iso & \RDerived\Gamma_X\left(\GtiltingObject^\vee \Ltimes \GtiltingObject \Ltimes \Sym^\bullet \Hom(V,S)^\vee \right). \eeqa Note that in fact the tensor products of these locally free sheaves do not need to be derived. It suffices to show then that the following bundle has no higher cohomology: $$B_{\alpha,\alpha'} := \Sigma^{-\alpha} S^\vee \otimes \underbrace{\Sigma^{\alpha'} S^\vee \otimes \Sym^\bullet \Hom(V,S)^\vee}_{P_{\alpha'}}.$$ We consider in particular the highlighted bundle $P_{\alpha'}$ with $$P_{\alpha'} := \Sigma^{\alpha'} S^\vee \otimes \Sym^\bullet \Hom(V,S)^\vee.$$ This bundle $P_{\alpha'}$ may be decomposed into terms $\Sigma^\mu S^\vee$ with positive weight $\mu \geq 0$: this follows immediately from the \emphasis{Littlewood-Richardson rule} for calculating tensor products of Schur powers \cite[Formula A.8]{Fulton:1996tk} using that $\alpha' \geq 0$. Similarly the whole bundle $B_{\alpha,\alpha'}$ may be decomposed into Schur powers $\Sigma^\mu S^\vee$ with $\mu \geq -\alpha \geq -\alpha_{top}$. The higher cohomology of these bundles vanishes by the proof of Proposition \ref{proposition.kapranov_exceptional_collection}, see \cite[Lemma 3.2(a)]{Kapranov:2009wf} for details.
\\
\step{2} We now show that $\XtiltingObject$ spans $D^b(X)$. By adjunction we have $$\Hom_X(\XtiltingObject, -) \iso \Hom_X(\bigBundleProjection^* \GtiltingObject, -) \iso \Hom_\Grassmannian(\GtiltingObject, \bigBundleProjection_* -).$$ Now $\bigBundleProjection$ is affine hence $\bigBundleProjection_*$ is injective. $\GtiltingObject$ is a tilting generator for $\Grassmannian$ hence by  Proposition \ref{proposition.tilting_generator_criterion} we have $\GtiltingObject^{\perp} \iso 0$. We deduce that $\XtiltingObject^{\perp} \iso 0$ as required.}
\end{proposition}

\begin{remark} Note that this construction does \remarkEmphasis{not} work with $X$ replaced by for instance the cotangent bundle $\cotangentBundle{\Grassmannian}$. In particular, for the simplest non-degenerate Grassmannian $\Grassmannian(2,4)$, it is noted in \cite[Remark 3.6(1)]{Kawamata:2005} that $$\Ext^2_{\cotangentBundle{\Grassmannian}}(\bigBundleProjection^* \GtiltingObject, \bigBundleProjection^* \GtiltingObject) \neq 0.$$ (We reuse the notation $\bigBundleProjection$ for the projection $\cotangentBundle{\Grassmannian} \to \Grassmannian$.)

For a construction of a tilting generator on $\cotangentBundle{\Grassmannian(2,4)}$ by another method see \cite{Toda:2008wl}.
\end{remark}

\subsection{Explicit descriptions}
\label{section.explicit_tilting_object}

\begin{samepage}
We restrict now to the $r=2$ case. The summands of $\mcT$ are as follows:
\begin{center} \begin{tikzpicture}
\bigWindowNodesSchurPowers{1.3}
\end{tikzpicture} \end{center}
\end{samepage}
\begin{samepage}
Expanding the Schur powers as in Example \ref{example.Schur_expansion_for_Grassmannian} we obtain:
\begin{center} \begin{tikzpicture}
\bigWindowNodes{1.3}
\end{tikzpicture} \end{center}
\end{samepage}

\section{Calculations on the tilting generator}
\label{section.calculations_on_generator}

\subsection{Outline}

We investigate now what the adjoint functors $L$ and $\rightAdjoint$ do to the summands of our tilting generator $\XtiltingObject$ from Appendix \ref{section.tilting_object} for the case $r=2$.

These results allow us to characterise the subcategory $\Im(L) \subset D^b(\twistBase)$ in Proposition \ref{lemma.generators_for_subcategory}. We also use them while understanding the action of the cotwist $C_F$ on $\Im L$ in Lemma \ref{proposition.image_within_window}, to apply $R$ to the convolution expressions for the $Fl^{\vee k}$ arising in Lemma \ref{proposition.generalized_Koszul_complexes}. All terms of each convolution, except the left-most term, are isomorphic to direct sums of the summands of $\XtiltingObject$. The left-most terms are isomorphic to directs sums of summands of $\XtiltingObject \otimes \O(1)$, and so we calculate the action of the adjoints on these too.

\begin{figure}[H]
\label{figure.generator_summands}
\caption{Summands of tilting generators $\XtiltingObject$ and $\XtiltingObject \otimes \O(1)$} 
\begin{center} \begin{tikzpicture}
\twoBigWindowNodes{1.3}
\twoBigWindowLabelling
\end{tikzpicture} \end{center}
\end{figure}

\subsection{Calculation of images under $L$ and $R$}

Our calculation is routine, but quite elaborate. We first calculate $\leftAdjoint (\O(k))$ for sheaves $\O(k)$ in the bottom row of the diagram:

\begin{lemma}\label{lemma.bottom_row_vanishing} We have $$\leftAdjoint (\O(k)) \iso \left\{ \begin{array}{cc} \O & k=0 \\ 0 & 0<k<d-1 \\ l^{\vee d-2} \otimes \det V^\vee [\dim \pi] & k=d-1 \end{array} \right\}.$$

\Proof{Using our setup $$\xymatrix{\hat{\maxStratum} \ar[r]^f \ar@/_{1em}/[rr]_j \ar[d]_\pi & \maxStratum \xyhook[r]^i & X \\ \twistBase & }$$ the first equality follows directly from the description of $\leftAdjoint$ in Proposition \ref{proposition.right_adjoint_description}: \beqa \leftAdjoint (\O) & = & \RDerived\pi_* (\omega_\pi \otimes \LDerived\composition^* \O) [\dim \pi] \\ & \iso & \RDerived\pi_* \omega_\pi [\dim \pi] \\ & \iso & \O. \eeqa Now we saw in Proposition \ref{proposition.resolutions_and_flattening} that $\pi$ is a projective bundle $\PP(V/H)$. In the $r=2$ case, $H$ is a line bundle so we write $l:=H$ as before and obtain \beqa \omega_\pi & \iso & \O_\pi(-d+1) \otimes \det \left( \frac{V}{l} \right)^\vee \\ & \iso & \O_\pi(-d+1) \otimes l \otimes \det V^\vee. \eeqa We now pull back our sheaves $\O(k)$ from $X$ to $\hat{\maxStratum}$ and write them in terms of the tautological bundle $\O_\pi(1)$, \beqa \LDerived\composition^* \O(k) = \O(k) & = & (\wedge^2 S^\vee)^k \\ & \iso & (S/l)^{\vee k} \otimes l^{\vee k} \\ & = & \O_\pi(k) \otimes l^{\vee k}, \eeqa using  $\wedge^2 S^\vee \iso l^\vee \otimes (S/l)^\vee$ which follows from the short exact sequence: $$0 \To (S/l)^\vee \To S^\vee \To l^\vee \To 0.$$ We then obtain \beqa \leftAdjoint (\O(k)) & \iso & \RDerived\pi_*(\omega_\pi \otimes \O_\pi(k) \otimes l^{\vee k} )[\dim \pi] \\ & \iso & \RDerived\pi_*\left(\O_\pi(k-d+1) \otimes l^{\vee k-1} \otimes \det V^\vee\right) [\dim \pi] \\ & \iso & \RDerived\pi_*\left(\O_\pi(k-d+1)\right) \\ & & \qquad \otimes l^{\vee k-1} \otimes \det V^\vee [\dim \pi], \eeqa and apply standard vanishing to give the result.
}
\end{lemma}

We have now applied $\leftAdjoint$ to the bottom row of our diagram. We apply $\leftAdjoint$ to the rest by a simple induction argument (Proposition \ref{proposition.action_of_right_adjoint}) which relies on the following lemma:

\begin{lemma}\label{lemma.induction_SES} We have $$0 \To \Sym^b S^\vee (a+1) \otimes l \To \Sym^{b+1} S^\vee (a) \To l^{\vee a+b+1} \otimes \O_\pi(a) \To 0$$ on $\hat{\maxStratum}$ where $a,b \geq 0$

\Proof{We have a short exact sequence $$0 \To (S/l)^\vee \To S^\vee \To l^\vee \To 0$$ which yields $$\xymatrix{0 \ar[r] & \Sym^b S^\vee \otimes (S/l)^\vee \ar[r] & \Sym^{b+1} S^\vee \ar[r] & l^{\vee (b +1)} \ar[r] & 0. \\ & \Sym^{b} S^\vee (1) \otimes l \xyequals[u] }$$ Multiplying this by $\O(a) \iso l^{\vee a} \otimes \O_\pi(a)$ we get the result.
}
\end{lemma}

We now perform our induction, according to the strategy shown below:

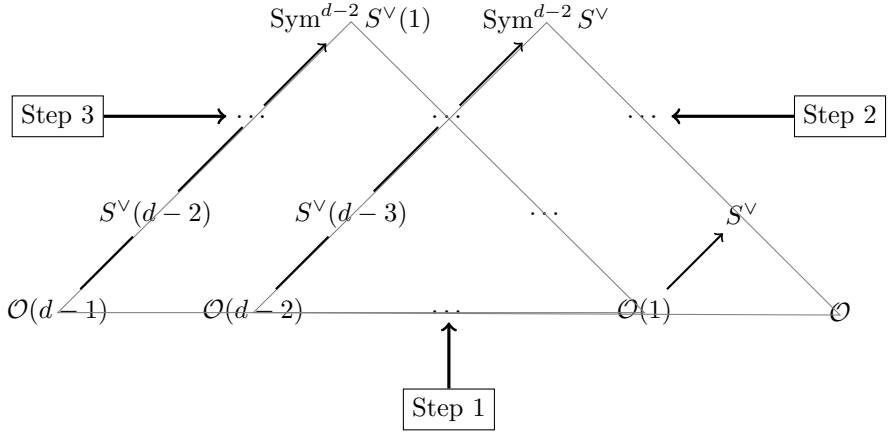
\begin{center} \begin{tikzpicture}
\twoBigWindowNodesArrows{1.3}
\inductionStrategyLabelling
\end{tikzpicture} \end{center}

\begin{proposition}\label{proposition.action_of_left_adjoint} For $0 \leq b \leq d-2$ we have $$\leftAdjoint (\Sym^b S^\vee (a)) \iso \left\{ \begin{array}{cc} 0 & a>0, \: \: a+b \leq d-2 \\ l^{\vee b} & a=0 \\ l^{\vee d-2-b} \otimes \det V^\vee [\dim \pi] & a+b=d-1 \end{array} \right\}.$$

\Proof{

\stepDescribed{1}{vanishing} Fix $a+b=k \leq d-2$. We prove the vanishing results first by increasing induction on $b$. Vanishing is known for $b=0$ by Lemma \ref{lemma.bottom_row_vanishing}. Now Lemma \ref{lemma.induction_SES} gives $$0 \To \LDerived\composition^* (\Sym^b S^\vee(a+1)) \otimes l \To \LDerived\composition^* (\Sym^{b+1} S^\vee(a)) \To l^{\vee a+b+1} \otimes \O_\pi(a) \To 0.$$ If $b < k$ then $0 < a \leq d-2$ and hence $\RDerived\pi_* (\O_\pi(a) \otimes \omega_\pi) \iso 0$  and we deduce $$\leftAdjoint (\Sym^b S^\vee(a+1)) \otimes l \iso \leftAdjoint (\Sym^{b+1} S^\vee(a))$$ which yields the vanishing by induction.

\stepDescribed{2}{right-most sheaves} Now for $b>0$, Lemma \ref{lemma.induction_SES} gives $$0 \To \LDerived\composition^* (\Sym^{b-1} S^\vee (1)) \otimes l \To \LDerived\composition^* (\Sym^b S^\vee) \To l^{\vee b} \To 0,$$ and now the vanishing just proved gives $$\leftAdjoint(\Sym^b S^\vee) \iso \pi_* (l^{\vee b} \otimes \omega_\pi) [\dim \pi] \iso l^{\vee b}$$ by the projection formula, as required.

\stepDescribed{3}{left-most sheaves} Here we prove the result by increasing induction on $b$. It holds for $b=0$. As before the proposition gives $$0 \To \LDerived\composition^* (\Sym^b S^\vee (d-1-b)) \otimes l \To \LDerived\composition^*(\Sym^{b+1} S^\vee (d-2-b))\To$$$$\qquad \qquad \To l^{\vee d-1} \otimes \O_\pi(d-2-b) \To 0.$$ This gives $$\leftAdjoint(\Sym^b S^\vee (d-1-b)) \otimes l \iso \leftAdjoint(\Sym^{b+1}  S^\vee (d-2-b)),$$ which allows us to complete the induction.
}
\end{proposition}

\begin{corollary}\label{proposition.action_of_right_adjoint} For $0 \leq b \leq d-2$ we have $$\rightAdjoint (\Sym^b S^\vee (a)) = \left\{ \begin{array}{cc} 0 & a>0, \: \: a+b \leq d-2 \\ l^{\vee b}[\dim \composition-\dim \pi] & a=0 \\ l^{\vee d-2-b} \otimes \det V^\vee [\dim \composition] & a+b=d-1 \end{array} \right\}.$$
\Proof{Use $R \iso L[\dim j - \dim \pi].$}
\end{corollary}

\begin{samepage}
Summarizing our results we have:
\begin{figure}[H]
\caption{Images under adjoints $\leftAdjoint$ and $\rightAdjoint$ of sheaves in Figure \ref{figure.generator_summands} (omitting the shifts and twists by $\det V$).} 
\begin{center} \begin{tikzpicture}
\twoBigWindowNodesResults{1.3}
\end{tikzpicture} \end{center}
\end{figure}
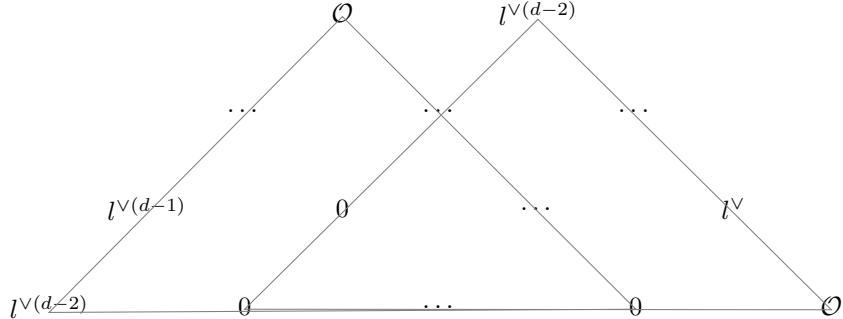
\end{samepage}

\subsection{Fullness of functor $L$}
\label{section.fullness}

\begin{lemma}\label{lemma.preimages_for_Homs_on_X}We have $$L(\Sym^b S^\vee) \iso l^{\vee b}$$ for $0 \leq b \leq d-2$ and the natural map $$\RHom_X(\Sym^a S^\vee, \Sym^b S^\vee) \overset{\phi_{a,b}}\To \RHom_{\twistBase}(l^{\vee a} , l^{\vee b})$$ induced by the functoriality of $L$ is surjective for $0 \leq a,b \leq d-2$.

\Proof{The first part comes from Proposition \ref{proposition.action_of_left_adjoint}. We now analyse the map $\phi_{a,b}$ to show surjectivity, as follows:

\step{1} Working as in Proposition \ref{proposition.base_tilting_object} we find that \beqa
\RHom_{\twistBase}(l^{\vee a} , l^{\vee b}) & \iso & \bigoplus_k \R\Gamma_{\PP V}(l^a \otimes l^{\vee b} \otimes \Sym^k(V \otimes l^\vee)) \\ & \iso & \bigoplus_k \R\Gamma_{\PP V}(l^{\vee b+k-a} \otimes \Sym^k V)
\eeqa
Similarly working as in Proposition \ref{proposition.tilting_object}, Step 1, we see that
\beqa \RHom_X(\Sym^a S^\vee, \Sym^b S^\vee) & \iso & \bigoplus_k \R\Gamma_{\Grassmannian}(\underbrace{\Sym^a S \otimes \Sym^b S^\vee \otimes \Sym^k (V \otimes S^\vee)}_{P_{a,b,k}}) \eeqa where $$P_{a,b,k} := \Sym^a S \otimes \Sym^b S^\vee \otimes \Sym^k (V \otimes S^\vee).$$ The map $\phi_{a,b}$ respects the summation and so it suffices to show that its $k^{\text{th}}$ summand, say $\phi_{a,b,k}$, is surjective.

\step{2} We identify a particular direct summand in $P_{a,b,k}$: the map $\phi_{a,b,k}$ will factor through this summand (after taking sections). First note that the $\Sym^k(\ldots)$ factor of $P_{a,b,k}$ can be decomposed into a direct sum of irreducibles by the Cauchy formula \cite[Theorem 2.3.2]{Weyman:wu}: we will only need that \beq{equation.Sym^k_decomposition} \Sym^k V \otimes \Sym^k S^\vee \into \Sym^k(V \otimes S^\vee).\eeq Similarly the Littlewood-Richardson rule \cite[Formula A.8]{Fulton:1996tk} can be used to decompose the complementary $\Sym^a S \otimes \Sym^b S^\vee$ factor of $P_{a,b,k}$: we just note that \beq{equation.Sym^a,b_decomposition} \Sym^{b-a} S^\vee \into \Sym^a S \otimes \Sym^b S^\vee.\eeq Observe that \eqref{equation.Sym^a,b_decomposition} does indeed make sense for $b-a < 0$, providing that we take $$\Sym^{-k} S^\vee := \Sym^k S.$$ Now putting \eqref{equation.Sym^k_decomposition} and \eqref{equation.Sym^a,b_decomposition} together, and using the Littlewood-Richardson rule once again, we have an inclusion $i_{a,b,k}$ as follows $$\xymatrix@C=4em{ \Sym^{b+k-a} S^\vee \otimes \Sym^k V \: \xyhook@<0.5ex>[r]^{\qquad \qquad i_{a,b,k}} & P_{a,b,k} \ar@<0.5ex>[l]^{\qquad \qquad \pi_{a,b,k}}},$$ as well as a corresponding projection, which we denote $\pi_{a,b,k}$.

\step{3} We can now describe the map $\phi_{a,b,k}$: it corresponds to $\R\Gamma_{\Grassmannian} (\pi_{a,b,k})$ under the following chain of isomorphisms: $$\R\Gamma_{\PP V}(l^{\vee b+k-a} \otimes \Sym^k V) \iso \Sym^{b+k-a} V^\vee \otimes \Sym^k V \iso \R\Gamma_{\Grassmannian} ( \Sym^{b+k-a} S^\vee \otimes \Sym^k V )  $$
Now $\R\Gamma_{\Grassmannian} (\pi_{a,b,k})$ is clearly a surjection, split by $\R\Gamma_{\Grassmannian} (i_{a,b,k})$, and hence the result follows.}
\end{lemma}


\bibliographystyle{\BibliographyLocation/Styles/halpha}
\bibliography{\BibliographyLocation/Bibliography.mod.stripped}

\end{document}